\documentclass{article}

\usepackage{natbib}

\usepackage{amsmath}
\usepackage{amsfonts}
\usepackage{amsthm}
\usepackage{graphicx}
\usepackage[latin1]{inputenc}
\usepackage{color}

\newtheorem{remark}{Remark}[section]
\newtheorem{theorem}[remark] {Theorem}

\newtheorem{lemma}[remark]{Lemma}
\newtheorem{cor}[remark]{Corollary}

\newtheorem{example}[remark]{Example}

\def\CF { \mathcal{F}}

\def\CP { \mathcal{P}}

\def\CN {\mathcal{N}}
\def\CS {\mathcal{S}}
\def\CG {\mathcal{G}}
\def\CP {\mathcal{P}}

\def\CA {\mathcal{A}}
\def\CE {\mathcal{E}}

\def\a {\alpha}

\def\th {\theta}

\def\a {\alpha}

\def\N {\mathbb{N}}

\def\N {\mathbb{N}}

\newcommand{\R}{\mathbb R}

\begin{document}

\title{Conditionally identically distributed species sampling
  sequences} 

\author{Federico Bassetti\footnote{Department of
    Mathematics, University of Pavia, via Ferrata 1, 27100 Pavia,
    Italy.\hspace{4,5cm} federico.bassetti@unipv.it}, Irene
  Crimaldi\footnote{Department of Mathematics, University of Bologna,
    Piazza di Porta San Donato 5, 40126 Bologna, Italy. 
    crimaldi@dm.unibo.it}, Fabrizio Leisen\footnote{Department of
    Economics, University of Pavia, Via San Felice, 5, 27100
    Pavia, Italy.\hspace{5cm} fabrizio.leisen@unimore.it} }

\date{}    

\maketitle
\abstract{Conditional identity in distribution (Berti et al. (2004))
  is a new type of dependence for random variables, which generalizes
  the well-known notion of exchangeability. In this paper, a class of
  random sequences, called {\em Generalized Species Sampling
    Sequences}, is defined and a condition to have conditional
  identity in distribution is given.  Moreover, a class of generalized
  species sampling sequences that are conditionally identically
  distributed is introduced and studied: the {\em Generalized Ottawa
    sequences} (GOS). This class contains a ``randomly reinforced''
  version of the P\'olya urn and of the Blackwell-MacQueen urn scheme.
  For the empirical means and the predictive means of a GOS, we prove
  two convergence results toward suitable mixtures of Gaussian
  distributions. The first one is in the sense of {\em stable}
  convergence and the second one in the sense of {\em almost sure
    conditional} convergence.  In the last part of the paper we study
  the length of the partition induced by a GOS at time $n$, i.e. the
  random number of distinct values of a GOS until time $n$. Under
  suitable conditions, we prove a strong law of large numbers and a
  central limit theorem in the sense of stable convergence. All the
  given results in the paper are accompanied by some examples.
  \\[5pt]
    }

\section{Introduction}\label{Sec.1} 

A sequence $(X_n)_{n\geq 1}$ of random variables defined on a
probability space $(\Omega,{\cal A},P)$ taking values in a Polish
space, is said a {\it species sampling sequence} if (a version) of the
 regular conditional distribution
of $X_{n+1}$ given $X(n):=(X_1,\dots,X_n)$ is the transition kernel
\begin{equation}\label{hanspit1}
K_{n+1}(\omega,\cdot):=
\textstyle\sum_{k=1}^{n} \tilde p_{n,k}(\omega)\delta_{X_k(\omega)}(\cdot)
+\tilde r_n(\omega)\mu(\cdot)
\end{equation}
where $\tilde p_{n,k}(\cdot)$ and $\tilde r_n(\cdot)$ are real--valued
measurable functions of $X(n)$ and $\mu$ is a probability measure. See
\cite{pitman2}.

As explained in \cite{han-pitman}, a species sampling
sequence $(X_n)_{n\geq 1}$ can be interpreted as the sequential random
sampling of individuals' species from a possibly infinite population
of individuals belonging to several species. If, for the sake of
simplicity, we assume that { $\mu$ is diffuse, then } the
interpretation is the following.  The species of the first individual
to be observed is assigned a random tag $X_1$, distributed according
to $\mu$. Given the tags $X_1, \dots X_n$ of the first $n$ individuals
observed, the species of the $(n+1)$-th individual is a new species
with probability $\tilde r_n$ and it is equal to the observed species
$X_k$ with probability $\sum_{j=1}^n \tilde p_{n,j}I_{\{X_j=X_k\}}$.

The concept of species sampling sequence is naturally related to that
of random partition induced by a sequence of observations.  Given a
random vector $X(n)=(X_1,\dots,X_n)$, we denote by $L_n$ the (random)
number of distinct values of $X(n)$ and by
$X^*(n)=(X^*_1,\dots,X^*_{L_n})$ the random vector of the distinct
values of $X(n)$ in the order in which they appear.  The {\em random
  partition induced by $X(n)$} is the random partition of the set 
$\{1,\dots,n\}$ given by \(
\pi^{(n)}=[\pi^{(n)}_1,\dots,\pi^{(n)}_{L_n}] \) where
$$
i \in \pi^{(n)}_k \Leftrightarrow X_i=X^*_k.
$$
{ Two distinct indices $i$ and $j$ clearly belong} to the same
block $\pi^{(n)}_k$ for a suitable $k$ if and only if $X_i=X_j$. It
follows that the {\em prediction rule} (\ref{hanspit1}) can be
rewritten as
\begin{equation}\label{hanspit2}
K_{n+1}(\omega,\cdot)=
{\textstyle\sum_{k=1}^{L_n(\omega)}} 
\tilde p_{n,k}^*(\omega) \delta_{X_k^*(\omega)}(\cdot)
+\tilde r_n(\omega)\mu(\cdot)
\end{equation}
where 
\[
\tilde p_{n,k}^*:=
{\textstyle\sum_{j \in \pi_k^{(n)}}} 
\tilde p_{n,j}.
\] 

In \cite{han-pitman} it is proved that if $\mu$ is diffuse and
$(X_n)_{n\geq 1}$ is an exchangeable sequence, then the coefficients
$\tilde p^*_{n,k}$ are almost surely equal to some function of
$\pi^{(n)}$ and they must satisfy a suitable recurrence relation.
Although there are only a few explicit prediction rules which give
rise to exchangeable sequences, this kind of prediction rules are
appealing for many reasons.  Indeed, exchangeability is a very natural
assumption in many statistical problems, in particular from the
Bayesian viewpoint, as well for many stochastic models.  Moreover,
remarkable results are known for exchangeable sequences: among others,
such sequences satisfy a strong law of large numbers and they can be
completely characterized by the well--known de Finetti representation
theorem. See, e.g., \cite{Aldous1985}.  Further, for an exchangeable
sequence the empirical mean $\sum_{k=1}^nf(X_k)/n$ and the predictive
mean, i.e.  $E[f(X_{n+1})|X_1,\dots,X_n]$, converge to the same limit
as the number of observations goes to infinity.  This fact can be
invoked to justify the use of the empirical mean in the place of the
predictive mean, which is usually harder to compute. Nevertheless, in
some situations the assumption of exchangeability can be too
restrictive.  For instance, instead of a classical P\'olya urn scheme,
it may be useful to deal with the so called randomly reinforced
P\'olya urn scheme. See, for example, \cite{cri}, \cite{cri-lei},
\cite{flournoy-may} and \cite{may-paganoni-secchi}.  Such a process
fails to be exchangeable but it can be still described with a
prediction rule which is not too far from (\ref{hanspit1}), see
Example \ref{ex:genpolya} of the present paper. Our purpose is 
to introduce and study a class of {\em generalized species sampling
  sequences}, which are generally not exchangeable but which still
have interesting mathematical properties.  

We thus need to recall the
notion of {\em conditional identity in distribution}, introduced and
studied in \cite{be-pra-ri}.  Such form of dependence generalizes the
notion of exchangeability preserving some of its nice predictive
properties.
One says that a sequence $(X_n)_{n\geq 1}$, defined on
$(\Omega,\CA,P)$ and taking values in a measurable space $(E,\CE)$, is
{\em conditionally identically distributed} with respect to a
filtration $\CG=(\CG_n)_{n\geq 0}$ (in the sequel, $\CG$-CID for
short), whenever $(X_n)_{n\geq 1}$ is $\cal G$--adapted and, for each
$n \geq 0$, $j \geq 1$ and every measurable real--valued bounded
function $f$ on $E$,
\[
E[f(X_{n+j})\,|\,\CG_n]=E[f(X_{n+1})\,|\,\CG_n]. 
\]
This means that, for each $n \geq 0$, all the random variables
$X_{n+j}$, with $j \geq 1$, are identically distributed conditionally
on $\CG_n$.  It is clear that every exchangeable sequence is a CID
sequence with respect to its natural filtration but a CID sequence is
not necessarily exchangeable.  Moreover, it is possible to show that a
$\cal G$--adapted sequence $(X_n)_{n \geq 1}$ is $\CG$-CID if and only
if, for each measurable real--valued bounded function $f$ on $E$,
\[
V_n^f:={\rm E}[f(X_{n+1})\,|\,{\CG}_n]  
\]
is a $\CG$--martingale, see \cite{be-pra-ri}. Hence, the sequence
$(V^f_n)_{n\geq 0}$ converges almost surely and in $L^1$ to a random
variable $V_f$. One of the most important features of CID sequences is
the fact that this random variable $V_f$ is also the almost sure
  limit of the empirical means. More precisely, CID sequences satisfy
  the following strong law of large numbers: for each real--valued bounded
  measurable function $f$ on $E$, the sequence $(M^f_n)_{n\geq 1}$, defined by
\begin{equation}\label{emp-mean}
M^f_n:=\frac{1}{n}\textstyle\sum_{k=1}^n f(X_k),
\end{equation} 
converges almost surely and in $L^1$ to $V_f$. It follows that also
the predictive mean ${\rm E}[f(X_{n+1})|X_1,\dots,X_n]$ converges
almost surely and in $L^1$ to $V_f$. In other words, CID sequences
share with exchangeable sequences the remarkable fact that the
predictive mean and the empirical mean merge when the number of
observations diverges. Unfortunately, while, for an exchangeable
sequence, we have $V_f=E[f(X_1)|{\cal T}]=\int f(x)m(\omega, {\rm
  d}x)$, where $\cal T$ is the tail--$\sigma$--field and $m$ is the
random directing measure of the sequence, it is difficult to
characterize explicitly the limit random variable $V_f$ for a CID
sequence. Indeed no representation theorems are available for CID
sequences.  See, e.g., \cite{aletti-may-secchi}.

The paper is organized as follows. In Section \ref{Sec.2} we state our
definition of generalized species sampling sequence, we discuss some
examples and we give a condition under which a generalized species
sampling sequence is CID with respect to a suitable filtration $\cal
G$. In Sections \ref{Sec.3} and \ref{Sec.4} we deal with a particular
class of generalized species sampling sequences which are CID: the
{\em generalized Ottawa sequences} (GOS for short).  We prove that,
for a GOS, under suitable conditions, the sequence
$\sqrt{n}(M_n^f-V_n^f)$ converges in the sense of {\em stable}
convergence to a mixture of Gaussian distributions.  Moreover, we
show that, under suitable conditions, also $\sqrt{n}(V_n^f-V_f)$
converges in the sense of {\em almost sure conditional convergence} to
another mixture of Gaussian distributions.  Both types of convergences
are stronger than the convergence in distribution.  These results are
accompanied by two examples. In Section \ref{Sec.5} we study the
length $L_n$ of the random partition induced by a GOS at time $n$,
i.e. the random number of the distinct values assumed by a GOS until
time $n$. In particular, a strong law of large numbers and a stable
central limit theorem are presented.  This section is also enriched by
some examples. The paper closes by a section devoted to proofs and by
an appendix in which the reader can find some results used for the
proofs.

\section{Prediction rules which generate a CID sequence}\label{Sec.2}

The Blackwell--MacQueen urn scheme provides the most
famous example of exchangeable prediction rule, that is
\[
P\{ X_{n+1} \in \cdot \,\,| X_1,\dots,X_n \}=
{\textstyle\sum_{i=1}^n} \frac{1}{\th+n}
\delta_{X_i}(\cdot) +\frac{\th}{\th+n}\mu(\cdot)
\] 
where $\th$ is a strictly positive parameter and $\mu$ is a
probability measure, see, e.g., \cite{black-mac} and \cite{pitman2}.
This prediction rule determines an exchangeable sequence $(X_n)_{n
  \geq 1}$ whose directing random measure is a Dirichlet process with
parameter $\th \mu(\cdot)$, see \cite{Ferguson1973}.  According to
this prediction rule, if $\mu$ is diffuse, a new species is observed
with probability $\th/(\th+n)$ and an old species $X_j^*$ is observed
with probability proportional to the cardinality of $\pi^{(n)}_j$, a
sort of {\em preferential attachment principle}.  This rule has its
analogous in term of random partitions in the so--called {\it Chinese
  restaurant process}, see \cite{pitman} and the references therein.

A {\em randomly reinforced prediction rule} of the same kind could
work as follows:
\begin{equation}\label{genchineserest}
P\{ X_{n+1} \in \cdot \,\,|X_1,\dots,X_n,Y_1,\dots,Y_n \}=
{\textstyle\sum_{i=1}^n} 
\frac{ Y_i}{\th + \sum_{j=1}^n Y_j}\delta_{X_i}(\cdot) + 
 \frac{\th}{\th + \sum_{j=1}^n Y_j}\mu(\cdot)
\end{equation}
where $\mu$ is a probability measure and $(Y_n)_{n \geq 1}$ is a sequence
of independent positive random variables.  If $\mu$ is diffuse, then
we have the following interpretation: each individual has a random
positive weight $Y_i$ and, given the first $n$ tags
$X(n)=(X_1,\dots,X_n)$ together with the weights
$Y(n)=(Y_1,\dots,Y_n)$, it is supposed that the species of the next
individual is a new species with probability $\th/(\th + \sum_{j=1}^n
Y_j)$ and one of the species observed so far, say $X_l^*$, with
probability $\sum_{i \in \pi^{(n)}_l} Y_i/(\th+\sum_{j=1}^n Y_j)$.
Again a preferential attachment principle.  Note that,
in this case, instead of describing the law of $(X_{n})_{n \geq 1}$
with the sequence of the conditional distributions of $X_{n+1}$ given
$X(n)$, we have a latent process $(Y_n)_{n \geq 1}$ and we
characterize $(X_{n})_{n \geq 1}$ with the sequence of the conditional
distributions of $X_{n+1}$ given $(X(n),Y(n))$.

Now that we have given an idea, let us formalize what we mean by
{\em generalized species sampling sequence}. Let $(\Omega,{\cal A},P)$
be a probability space and $E$ and $S$ be two Polish spaces, endowed with
their Borel $\sigma$-fields $\CE$ and $\CS$, respectively. 
 In the sequel, $\CF^Z=(\CF^Z_n)_{n \geq 0}$ will stand for the
 natural filtration associated with any sequence of random variables 
$(Z_n)_{n \geq 1}$ on  $(\Omega,{\cal A}, P)$ and we set 
$\CF^Z_\infty=\vee_{n \geq 0} \CF^Z_n $. Finally, 
$\CP_n$ will denote the set of all partitions of $\{1,\dots,n\}$.

We shall say that a sequence $(X_n)_{n\geq 1}$ of random variables on
$(\Omega,{\cal A}, P)$, with values in $E$, is a { generalized
  species sampling sequence} if:
\begin{itemize}
\item {($h_1$) $X_1$ has distribution $\mu$.}
\item{($h_2$) There exists a sequence $(Y_n)_{n\geq 1}$ of random
    variables with values in $(S,\CS)$ such that, for each $n\geq 1$,
    a version of the regular conditional distribution of $X_{n+1}$ given
$$
{\cal F}_n:=
{\cal F}_n^X \vee{\cal F}_n^Y
$$
is 
\begin{equation}\label{cond-distr}
K_{n+1}(\omega,\cdot)=
\textstyle\sum_{i=1}^n p_{n,i}(\pi^{(n)}(\omega),Y(n)(\omega))
\delta_{X_i(\omega)}(\cdot) 
+r_n(\pi^{(n)}(\omega),Y(n)(\omega)) \mu(\cdot)
\end{equation}
with $p_{n,i}(\cdot,\cdot)$ and $r_n(\cdot,\cdot)$ suitable measurable 
functions defined on $\CP_n \times S^n$ with values in $[0,1]$. 
}
\item{($h_3$) $X_{n+1}$ and $(Y_{n+j})_{j\geq 1}$ are conditionally 
independent given $ {\cal F}_n$.}
\end{itemize}

\begin{example}\label{ex-2-gen}
  \rm Let $\mu$ be a probability measure on $E$, $(\nu_{n})_{n\geq 1}$
  be a sequence of probability measures on $S$, $(r_n)_{n\geq 1}$ and
  $(p_{n,i})_{n \geq 1, \, 1\leq i\leq n}$ be measurable functions
  such that 
$$r_n: \CP_n \times S^n \to [0,1],\qquad p_{n,i}:\CP_n \times
  Z^n \to [0,1]$$ 
and
\begin{equation}\label{condpart1}
{\textstyle\sum_{i=1}^n} 
p_{n,i}(q_n,y_1,\dots,y_n)+r_n(q_n,y_1,\dots,y_n)=1
\end{equation}
for each $n \geq 1$ and each $(q_n,y_1,\dots,y_n)$ in $\CP_n \times
S^n$. By the Ionescu Tulcea Theorem, there are two sequences of random
variables $(X_n)_{n\geq 1}$ and $(Y_{n})_{n\geq 1}$, defined on a
suitable probability space $(\Omega,{\cal A},P)$, taking values in $E$
and $S$ respectively, such that conditions $(h_1)$, $(h_2)$ and the
following condition are satisfied:
\begin{itemize}
\item{$Y_{n+1}$ has distribution $\nu_{n+1}$ and it is
independent of the $\sigma$-field
$$
{\CF}_n\vee\sigma(X_{n+1})=
{\CF}_{n+1}^X\vee{\CF}_n^Y.
$$
}
\end{itemize}
This last condition implies that, for each $n$, $(Y_{n+j})_{j\geq 1}$ is
independent of ${\cal F}_{n+1}^X\vee{\cal F}_n^Y$. It follows, in
particular, that $(Y_n)_{n\geq 1}$ is a sequence of independent random
variables. Therefore, also $(h_3)$ holds true. Indeed,
 for each real--valued bounded ${\cal F}_n$-measurable random
variable $V$, each bounded Borel function $f$ on $E$, each $j\geq
1$ and each bounded Borel function $h$ on $S^j$, we have
\begin{equation*}
\begin{split}
{\rm E}[Vf(X_{n+1})h(Y_{n+1},\dots,Y_{n+j})]&=
{\rm E}\big[\,Vf(X_{n+1})
{\rm E}[h(Y_{n+1},\dots,Y_{n+j})\,|\,{\cal F}_n\vee\sigma(X_{n+1})]\,
\big]\\
&=
{\rm E}[\,Vf(X_{n+1})
{\textstyle\int}h(y_{n+1},\dots,y_{n+j})\,
\nu_{n+1}({\rm d}y_{n+1})\dots({\rm d}y_{n+1})\,]\\
&=
{\rm E}\big[\,V{\rm E}[f(X_{n+1})\,|\,{\cal F}_n]
{\textstyle\int}h(y_{n+1},\dots,y_{n+j})\,
\nu_{n+1}({\rm d}y_{n+1})\dots({\rm d}y_{n+1})\,\big].
\end{split}
\end{equation*}
On the other hand, we have
$$
{\rm E}[h(Y_{n+1},\dots,Y_{n+j})\,|\,{\cal F}_n]
={\textstyle\int}h(y_{n+1},\dots,y_{n+j})\,
\nu_{n+1}({\rm d}y_{n+1})\dots({\rm d}y_{n+1})
$$
hence
$$
{\rm E}[f(X_{n+1})h(Y_{n+1},\dots,Y_{n+j})\,|\,{\cal F}_n]
=
{\rm E}[f(X_{n+1})\,|\,{\cal F}_n]
{\rm E}[h(Y_{n+1},\dots,Y_{n+j})\,|\,{\cal F}_n].
$$
This fact is sufficient in order to conclude that also
assumption $(h_3)$  is verified.$\qquad\Diamond$
\end{example}

In order to state our first result concerning generalized species
sampling sequences, we need some further notation.
Set
\[
\begin{split}
& p^{*}_{n,j}(\pi^{(n)})=p^{*}_{n,j}(\pi^{(n)},Y(n)):=
{\textstyle\sum_{i \in \pi^{(n)}_j}}p_{n,i}(\pi^{(n)},Y(n)) 
\quad \hbox{for } j=1,\dots,L_n 
\\
\mbox{and} \quad & \\
&r_n:=r_n(\pi^{(n)},Y(n)).  
\\
\end{split}
\]
Given a partition $\pi^{(n)}$, denote by $[\pi^{(n)}]_{j+}$ the
partition of $\{1,\dots,n+1\}$ obtained by adding the element $(n+1)$
to the $j$-th block of $\pi^{(n)}$. Finally, denote by
$[\pi^{(n)};(n+1)]$ the partition obtained by adding a block
containing $(n+1)$ to $\pi^{(n)}$.  For instance, if
$\pi^{(3)}=[(1,3);(2)]$, then $[\pi^{(3)}]_{2+}=[(1,3);(2,4)]$ and
$[\pi^{(3)};(4)]=[(1,3);(2);(4)]$.

\begin{theorem}\label{cond-cidness}
    A generalized species sampling sequence $(X_n)_{n \geq 1}$ 
    with $\mu$ diffuse is a CID sequence with respect to the
  filtration $\CG=(\CG_n)_{n \geq 0}$ with ${\CG}_n:={\cal F}_n^X \vee{\cal
    F}_{\infty}^Y $ if and only if, for each $n$, the following
  condition holds $P$-almost surely:
\begin{equation}\label{cond-cid}
p_{n,j}^*(\pi^{(n)})=
r_n p_{n+1,j}^*([\pi^{(n)};\{n+1\}])
+{\textstyle\sum_{l=1}^{L_n}} 
p_{n+1,j}^*([\pi^{(n)}]_{l+})p_{n,l}^*(\pi^{(n)})
\end{equation}
for $1\leq j \leq L_n$. 
\end{theorem} 

The next example generalizes the well-known two parameter
Poisson-Dirichlet process.

\begin{example}\label{ex-00} {\rm Let $\th>0$ and
    $\alpha\geq 0$. Moreover, let {  $\mu$ be a
      probability measure on $E$ and}, $(\nu_n)_{n \geq 1}$ be a
    sequence of probability measures on $(\a,+\infty)$.  Consider the
    following sequence of functions
\[
\begin{split}
p_{n,i}(q_n,y(n))&:=\frac{y_i-\a/C_{i}(q_n)}{\th + \sum_{j=1}^n y_j} \qquad
 \\
r_n(q_n,y(n))&:=\frac{\th+\a L(q_n)}{\th + \sum_{j=1}^n y_j}
\end{split}
\]
where $y(n)=(y_1,\dots,y_n)\in (\a,+\infty)^n, \,\, q_n \in \CP_n$,
$C_{i}(q_n)$ is the cardinality of the block in $q_n$ which contains
$i$ and $L(q_n)$ is the number of blocks of $q_n$.  It is easy to see 
that such functions satisfy (\ref{condpart1}). Hence, by Example
\ref{ex-2-gen}, there exists 
 a generalized species sampling sequence 
$(X_n)_{n \geq 1}$ for which
\begin{equation}\label{genchineserest2}
P\{ X_{n+1} \in \cdot \,\,|X(n),Y(n)\}=
{\textstyle \sum_{l=1}^{L_n}} 
\frac{ \sum_{i \in \pi_l^{(n)}} Y_i-\a }
{\th + \sum_{j=1}^n Y_j}\delta_{X_l^*}(\cdot) + 
\frac{\th +\a L_n }{\th + \sum_{j=1}^n Y_j}\mu(\cdot).
\end{equation}
where $(Y_n)_{n \geq 1}$ is a sequence of independent random variables such
that each $Y_n$ has law $\nu_n$.
{  If $\mu$ is diffuse,} one can easily check   
that (\ref{cond-cid}) of Theorem  
\ref{cond-cidness} holds and so $(X_n)_{n \geq 1}$ is a CID sequence 
with respect to $\CG=({\cal F}_n^X \vee{\cal F}_{\infty}^Y)_{n \geq 1}$. 

It is  worthwhile
noting that if $Y_n=1$ for every $n \geq 1$ and $\a$
 belongs to $[0,1]$,  then  we get an exchangeable sequence directed by the 
well-known  two parameter Poisson-Dirichlet process: i.e. an exchangeable 
sequence described by the prediction  rule
\[
P\{ X_{n+1} \in \cdot \,\,|X_1,\dots,X_n \}=
{\textstyle\sum_{l=1}^{L_n}} 
\frac{|\pi_l^{(n)}|-\a}{\th + n} \delta_{X_l^*}(\cdot) + 
 \frac{\th+\a L_n}{\th + n}\mu(\cdot).
\]
See, e.g., \cite{PitmanYor1997} and \cite{pitman}}.$\qquad\Diamond$
\end{example}

  A special case of the previous example is the randomly reinforced
  Blackwell-McQueen urn scheme (\ref{genchineserest}). However this
  prediction rule may be collocated in a more general class of
  generalized species sampling sequences, that are CID. In the next
  sections, we shall introduce and study this class, called
  ``Generalized Ottawa Sequences''.

\section{Generalized Ottawa sequences}\label{Sec.3}

We shall say that a generalized species sampling sequence 
$(X_n)_{n \geq 1}$
is  a {\em generalized Ottawa sequence} or, more briefly, a
{GOS}, if for every $n \geq 1$
\begin{itemize}
\item  The functions $r_n$ and $p_{n,i}$  ($i=1,\dots,n$)  do not
  depend on the partition, 
hence
\begin{equation}\label{cond-distr2}
K_{n+1}(\omega,\cdot)=
\textstyle\sum_{i=1}^n p_{n,i}(Y(n)(\omega))
\delta_{X_i(\omega)}(\cdot)
+r_n(Y(n)(\omega)) \mu(\cdot).
\end{equation}
\item  The functions $r_n$ are striclty positive and 
\begin{equation}\label{gos-functions}
r_n(Y_1,\dots,Y_n)\geq r_{n+1}(Y_1,\dots,Y_n,Y_{n+1})
\end{equation}
almost surely.
\item The functions $p_{n,i}$ satisfy
\begin{equation}\label{gos-functions2}
\begin{split}
p_{n,i}&:=\frac{r_n}{r_{n-1}}\,p_{n-1,i}\qquad\hbox{for } i=1,\dots,n-1\\
p_{n,n}&:=1-\frac{r_n}{r_{n-1}}
\end{split}
\end{equation}
with $r_0=1$.
\end{itemize}

For simplicity, from now on, we shall denote by $r_n$ and
$p_{n,i}$ the $\CF^Y_n$-measurable random variables $r_n(Y(n))$ and
$p_{n,i}(Y(n))$, that is
$r_n:=r_n(Y(n))$ and $ p_{n,i}:=p_{n,i}(Y(n))$.

  First of all let us stress that any GOS is a CID sequence with
  respect to the filtration $\CG=(\CF^X_n\vee\CF^Y_{\infty})_{n\geq
    0}$.  Indeed, since $\CG_n=\CF_n\vee\sigma(Y_{n+j}:j\geq 1)$,
  condition $(h3)$ implies that
\begin{equation}\label{legame-fondamentale}
{\rm E}[f(X_{n+1})\,|\,{\CG}_n]={\rm E}[f(X_{n+1})\,|\,{\cal F}_n]
\end{equation}
for each bounded Borel function $f$ on $E$
and hence, by $(h2)$, one gets
$$
V^f_n:={\rm E}[f(X_{n+1})\,|\,{\CG}_n]=
\textstyle\sum_{i=1}^n p_{n,i}f(X_i)+r_n {\rm E}[f(X_1)].
$$
 Since the random
variables $p_{n+1,i}$ are ${\CG}_n$-measurable it follows
that
\begin{equation*}
\begin{split}
{\rm E}[V^f_{n+1}\,|\,{\CG}_n]&=
{\textstyle\sum_{i=1}^{n}}p_{n+1,i}f(X_i)+
p_{n+1,n+1}{\rm E}[f(X_{n+1})\,|\,{\CG}_n]
+r_{n+1}{\rm E}[f(X_1)]\\
&=
\frac{r_{n+1}}{r_n}
{\textstyle\sum_{i=1}^{n}}p_{n,i}f(X_i)+
V_n^f-\frac{r_{n+1}}{r_n}V^f_n
+r_{n+1}{\rm E}[f(X_1)]\\&
=\frac{r_{n+1}}{r_n}V^f_n-
r_{n+1}{\rm E}[f(X_1)]
+
V^f_n-
\frac{r_{n+1}}{r_n}V^f_n
+r_{n+1}{\rm E}[f(X_1)]=V^f_n.
\end{split}
\end{equation*}

 Some examples follow.

\begin{example}\label{ex-2bis} 
\rm Consider a GOS for which 
\[
Y_n=a_n
\]
where $(a_n)_{n\geq 0}$ is a {\it decreasing numerical}
 sequence with $a_0=1$,  $a_n>0$ and $r_n(y_1,\dots,y_n)=y_n$.$\qquad\Diamond$
\end{example}

\begin{example}\label{ex-2tris} 
\rm  
Let $(Y_n)_{n\geq 1}$ be a Markov chain taking values in  $(0,1]$, with
$Y_1=1$ and transition probability kernel given by
\[
P\{Y_{n+1} \leq x|Y_{n}\}= 
\frac{x}{Y_n} I_{(0,Y_n)}(x)+I_{[Y_n,+\infty)}(x)
\quad n \geq 1.
\] 
Then we have $Y_n\geq Y_{n+1}$ a.s. for all $n\geq 1$. Thus we can
consider a GOS with $r_n(y_1,\dots,y_n)=y_n$.  $\Diamond$
\end{example}

As we shall see in the next example, the randomly reinforced
Blackwell--McQueen urn scheme gives rise to a GOS.

\begin{example}\label{ex-2} \rm 
  Let $\mu$ be a probability measure on $E$, $(\nu_{n})_{n\geq 1}$ be
  a sequence of probability measures on $S$ and $(r_n),\, (p_{n,i})$
  measurable functions as in (\ref{gos-functions}) and
  (\ref{gos-functions2}).  Following Example \ref{ex-2-gen}, 
there exist two sequences of random variables
  $(X_n)_{n\geq 1}$ and $(Y_{n})_{n\geq 1}$, defined on a suitable
  probability space $(\Omega,{\cal A},P)$,   such that
    each $Y_n$ has law $\nu_n$ and it is independent of
    ${\CF}^X_{n}\vee{\CF}^Y_{n-1}$ and $(X_n)_{n\geq 1}$ follows the
    prediction rule (\ref{cond-distr2}), i.e. it is a GOS.

    As special case one can consider $S=\R_+$ and
$$
r_n(y_1,\dots,y_n)=
\frac{\th}{\th+\sum_{j=1}^n y_j}
$$ 
with $\th >0$.$\qquad\Diamond$
\end{example} 

Particular case of the previous example is the following randomly
reinforced P\'olya urn.

\begin{example}[A randomly reinforced P\'olya urn]\label{ex:genpolya} {\rm An
    urn contains $b$ black and $r$ red balls, $b$ and $r$ being
    strictly positive integer numbers.  Repeatedly (at each time
    $n\geq 1$), one ball is drawn at random from the urn and then
    replaced together with a positive random number $Y_n$ of
    additional balls of the same color.  For each $n$, the random
    number $Y_n$ must be independent of the preceding numbers and of
    the drawings until time $n$.  If we denote by $X_n$ the indicator
    function of the event $\{\hbox{black ball at time } n\}$, then we
    clearly have $E=\{0,1\}$,
$$
\mu(0)=\frac{r}{b+r},\qquad \mu(1)=\frac{b}{b+r},\qquad 
$$
and 
$$
P\{ X_{n+1} \in \cdot \,\,|X(n),Y(n)\}
= 
\frac{1}{b+r+{\textstyle\sum_{j=1}^n} Y_j}
{\textstyle\sum_{i=1}^n}Y_i\delta_{X_i(\omega)}(\cdot)
+\frac{b+r}{b+r+{\textstyle\sum_{j=1}^n} Y_j}\mu(\cdot).
$$
Note that the sequence $(X_n)_{n\geq 1}$ is generally not
exchangeable. Indeed, it is straightforward to prove that, even if the
random variables $Y_n$ are identically distributed, the sequence
$(X_n)_{n\geq 1}$ is not exchangeable (apart from particular cases).
}$\qquad\Diamond$
\end{example}

\section{Convergence results for a GOS}\label{Sec.4}

In this section we prove some limit theorems for a GOS under stable
convergence and almost sure conditional convergence.

Stable convergence has been introduced by \cite{re} and subsequently
studied by various authors,  see, for example, \cite{aldous-eagleson},
  \cite{jacod-memin},  \cite{hall-heyde}. A detailed treatment,
including some strengthened forms of stable convergence, can be found
in \cite{cri-le-pra}.

Given a probability space $(\Omega,{\cal
  A}, P)$ and a Polish space $E$ (endowed with its Borel
$\sigma$-field ${\cal E}$), a kernel $K$ on $E$ is a family
$K=(K(\omega,\cdot))_{\omega\in\Omega}$ of probability measure on $E$
such that, for each bounded Borel function $g$ on $E$, the function
$$
K(g)(\omega)=\textstyle\int g(x) K(\omega,{\rm d}x)
$$  
is measurable with respect to $\cal A$. Given a sub-$\sigma$-field
$\cal H$ of $\cal A$, we say that the kernel $K$ is $\cal
H$-measurable if, for each bounded Borel function $g$ on $E$, the
random variable $K(g)$ is measurable with respect to $\cal H$. In the
following, the symbol $\cal N$ will denote the sub-$\sigma$-field
generated by the $P$-negligible events of $\cal A$. Given a
sub-$\sigma$-field $\cal H$ of $\cal A$ and a ${\cal H}\vee{\cal
  N}$-measurable kernel $K$ on $E$, a sequence $(Z_n)_{n \geq 1}$ of
random variables on $(\Omega,{\cal A},P)$ with values in $E$ converges
{\it$\cal H$-stably} to $K$ if, for each bounded continuous function
$g$ on $E$ and for each ${\cal H}$--measurable real--valued bounded
random variable $W$
\begin{equation*}
{\rm E}[g(Z_n)\, W]\longrightarrow{\rm E}[K(g)\, W].
\end{equation*}
 If $(Z_n)_{n \geq 1}$ converges $\cal H$-stably to $K$
  then, for each $A\in{\cal H}$ with $P(A)\neq 0$,  
  the sequence $(Z_n)_{n \geq 1}$ converges in distribution under the 
  probability measure $P_A=P(\cdot|A)$ to the probability measure  $P_AK$ on
  $E$ given by
\begin{equation}\label{mistura2}
P_AK(B)=P(A)^{-1}{\rm E}[I_A K(\cdot,B)]=
\textstyle\int K(\omega,B)\,P_A({\rm d}\omega)
 \qquad\hbox{for each } B \in \CE.
\end{equation}
In particular, if $(Z_n)_{n \geq 1}$ converges $\cal H$-stably to $K$,
then $(Z_n)_{n \geq 1}$ converges in distribution to the probability measure  
$PK$ on $E$ given by
\begin{equation}\label{mistura}
PK(B)={\rm E}[K(\cdot,B)]=\textstyle\int K(\omega,B)\,P({\rm d}\omega)
 \qquad\hbox{for each } B \in \CE.
\end{equation}
Moreover, if all the random variables $Z_n$ are $\cal H$-measurable,
then the $\cal H$-stable convergence obviously implies the $\cal A$-stable
convergence. 

Given a filtration $\CG=(\CG_n)_{n\geq 0}$ and a kernel $K$ on $E$, we
shall say that, with respect to $\CG$, the sequence $(Z_n)_{n\geq 1}$
converges to $K$ in the sense of the almost sure conditional convergence
if, for each   bounded continuous function $g$, we have
\begin{equation*}
{\rm E}[g(Z_n)\,|\, {\CG}_n]\longrightarrow K(g)
\qquad\hbox{almost surely.}
\end{equation*}
If $(Z_n)_{n\geq 1}$ converges to $K$ in the
sense of the almost sure conditional convergence with respect to a
filtration $\CG$, then $(Z_n)_{n\geq 1}$ also converges
${\CG}_{\infty}$-stably to $K$, see \cite{cri}.

Throughout the paper, if $U$ is a positive random variable, we shall call the
Gaussian kernel associated with $U$ the family
$$
{\cal N}(0,U)=\big({\cal N}(0,U(\omega))\big)_{\omega\in\Omega}
$$
of Gaussian distributions with zero mean and variance equal to
$U(\omega)$ (with ${\cal N}(0,0):=\delta_0$). Note that, in this case,
the probability measure defined in (\ref{mistura2}) and
(\ref{mistura}) is a mixture of Gaussian distributions.

It is worthwhile to recall that, if
  $(X_n)_{n \geq 1}$ is a GOS, then it is a CID sequence with respect
  to the filtration $\CG=(\CF^X_n\vee\CF^Y_{\infty})_{n\geq 0}$ (as
  shown in Section \ref{Sec.3}) and so the sequence $V_n^f$ (defined
  in section \ref{Sec.3}) converges almost surely and in $L^1$ to a
  random variable $V_f$, whenever $f$ is a bounded Borel
  function on $E$. Moreover, the random variable $V_f$ is also the
  almost sure (and in $L^1$) limit of the empirical mean
$$
M_n^f=\frac{1}{n}{\textstyle\sum_{k=1}^n}f(X_k).
$$
We are ready to state the main theorems of this section.

\begin{theorem}\label{clt2} {  Let $(X_n)_{n\geq 1}$ be a GOS. 
Using the above notation, for each bounded Borel
      function $f$ and each $n\geq 1$, let us set
\begin{equation*}
S^f_n=\sqrt{n}(M^f_{n}- V^f_n)
\end{equation*}
and, for $1\leq j\leq n$,
\begin{equation*}
Z^f_{n,j}=
\frac{1}{\sqrt{n}}\big[f(X_j)-jV^f_j+(j-1)V^f_{j-1}\big]
=\frac{1}{\sqrt{n}}(1+jp_{j,j})\big[f(X_j)-V^f_{j-1}\big].
\end{equation*}
Suppose that:
\\[5pt]
\indent (a) $U^f_n:=\sum_{j=1}^n(Z^f_{n,j})^2\stackrel{P}\longrightarrow U_f$.
\\[5pt]
\indent (b) $(Z^f_n)^*:=
\sup_{1\leq j\leq n}|Z^f_{n,j}|\stackrel{L^1}\longrightarrow 0$.
\\[5pt]
\noindent Then the sequence $(S^f_{n})_{n \geq 1}$ converges $\cal A$-stably
to the Gaussian kernel ${\cal N}(0,U_f)$.
\\[5pt]
\indent In particular, condition (a) and (b) are satisfied if the
following conditions hold:
\\[5pt]
\indent (a1) $U^f_n\stackrel{a.s.}\longrightarrow U_f$.
\\[5pt]
\indent (b1) $\sup_{n \geq 1}{\rm E}[(S^f_n)^2]<+\infty$.  }
\end{theorem}

Let us see an application of the previous theorem 
in the next example.

\begin{example}\label{example-clt2} 
{\rm Let us consider the setting  of Example \ref{ex-2} with 
$$
r_k=\frac{\th}{\th+\textstyle\sum_{i=1}^k Y_i}.  
$$ 
where $\th>0$ and the random variables $Y_n$ are identically
distributed with $Y_n\geq\gamma>0$ and ${\rm E}[Y_n^4]<+\infty$.
Given a bounded Borel function $f$ on $E$ we are going to prove that
the sequence $(S^f_n)_{n\geq 1}$ (defined in Theorem \ref{clt2})
converges $\cal A$-stably to the Gaussian kernel
$${\cal N}\big(0,\Delta
(V_{f^2}-V_f^2)\big),$$ 
where $\Delta:={\rm Var}[Y_1]/{\rm E}[Y_1]^2$.

Without loss of generality, we may assume that $f$ takes values in $[0,1]$. 
Let us observe that, after some calculations, we have
\begin{equation*}
S_n^f=\frac{1}{\sqrt{n}}\left({\textstyle\sum_{i=1}^n}f(X_i)-n\,
\frac{\th{\rm E}[f(X_1)]+{\textstyle\sum_{i=1}^n}Y_if(X_i)}
{\th+{\textstyle\sum_{j=1}^n}Y_j}\right).
\end{equation*}
If we set $b:=\th{\rm E}[f(X_1)]$ and 
$\widehat{Y_i}:=Y_i-{\rm E}[Y_i]=Y_i-m$,
then we can write
\begin{equation*}
S_n^f=\frac{1}{\sqrt{n}}\frac{1}{(\th+{\textstyle\sum_{j=1}^n}Y_j)}
\big[(\th+{\textstyle\sum_{j=1}^n}\widehat{Y_j}) {\textstyle\sum_{i=1}^n}f(X_i)
-nb - n{\textstyle\sum_{i=1}^n}\widehat{Y_i}f(X_i)\big].
\end{equation*}
Therefore, since $Y_n\geq \gamma$ and $0\leq f(X_n)\leq 1$ for each
$n$, we obtain
\begin{equation*}
\begin{split}
{\rm E}[(S_n^f)^2]&\leq\frac{2n}{(\th+\gamma n)^2}
\big({\rm E}[\,(\th+{\textstyle\sum_{j=1}^n}\widehat{Y_j})^2\,]+
{\rm E}[\,(b+{\textstyle\sum_{i=1}^n}\widehat{Y_i}f(X_i)\,)^2\,]\big)\\
&\leq 
\frac{2n}{(\th+\gamma n)^2}\big(\th^2+b^2+2n{\rm Var}[Y_1]\big)\leq C
\end{split}
\end{equation*}
where $C$ is a suitable constant. Finally, let us observe that, after some
calculations, we get
\begin{equation*}
\begin{split}
U_n^f&=\frac{1}{n}{\textstyle\sum_{j=1}^n}
\big[f(X_j)-jV^f_j+(j-1)V^f_{j-1}\big]^2\\
&=
\frac{1}{n}{\textstyle\sum_{j=1}^n}
\big[f(X_j)+A^f_jY_j-B_jY_jf(X_j)-V^f_{j-1}\big]^2,
\end{split}
\end{equation*}
where
\begin{equation*}
\begin{split}
A^f_j&=j
\big[(\th+{\textstyle\sum_{i=1}^j}Y_i)
(\th+{\textstyle\sum_{i=1}^{j-1}}Y_i)\big]^{-1}
\big[b-\th f(X_j)+{\textstyle\sum_{i=1}^{j-1}}Y_if(X_i)\big],\\
B_j&=j\big[(\th+{\textstyle\sum_{i=1}^j}Y_i)
(\theta+{\textstyle\sum_{i=1}^{j-1}}Y_i)\big]^{-1}
{\textstyle\sum_{i=1}^{j-1}}Y_i.
\end{split}
\end{equation*}
Hence, we have
\begin{equation*}
\begin{split}
&U^f_n=\frac{1}{n}{\textstyle\sum_{j=1}^n}
\big[f^2(X_j)+(A^f_j)^2Y_j^2+B_j^2Y_j^2f^2(X_j)+(V^f_{j-1})^2
-2f(X_j)V^f_{j-1}\big]
\\
&+\!\frac{2}{n}{\textstyle\sum_{j=1}^n}
\big[A^f_jY_jf(X_j)-B_jY_jf^2(X_j)-A^f_jY_jV^f_{j-1}
-A^f_jB_jY_j^2f(X_j)+B_jY_jf(X_j)V^f_{j-1}\big].
\end{split}
\end{equation*}
Recall that we have the following almost sure convergences: 
\begin{equation*}
f^q(X_n)/n\longrightarrow 0\,(\hbox{for }q=1,2),\quad
V^f_{n}\longrightarrow V_f
\end{equation*}
\begin{equation*}
\frac{1}{n}{\textstyle\sum_{j=1}^n}Y_j^{r}\longrightarrow 
{\rm E}[Y_1^r]\,(\hbox{for }r=1,2),\quad
\frac{1}{n}{\textstyle\sum_{j=1}^n} f(X_j)\longrightarrow V_f,\quad
\frac{1}{n}{\textstyle\sum_{j=1}^n} f^2(X_j)\longrightarrow V_{f^2}.
\end{equation*}
>From the above
relations, we get
\begin{equation*}
B_j\stackrel{a.s.}\longrightarrow 1/{\rm E}[Y_1].
\end{equation*}
In order to study the convergence of 
$\frac{1}{n}{\textstyle\sum_{j=1}^n}Y_j^{r}f^q(X_j)$ for $r,q=1,2$, 
let us set
$$
Z_n:={\textstyle\sum_{j=1}^n}
\frac{1}{j}\big(Y_j^{r}f^q(X_j)-{\rm E}[Y_j^{r}f^q(X_j)\,|\, 
{\cal F}_{j-1}
]\big).
$$
The sequence $(Z_n)_{n \geq 1}$ is a martingale with respect to 
${\cal F}=(\CF_n)_{n \geq 1}$ 
such that
\begin{equation*}
\begin{split}
{\rm E}[Z_n^2]&={\textstyle\sum_{j=1}^n}
\frac{1}{j^2}{\rm E}\big[\,
(Y_j^{r}f^q(X_j)-{\rm E}[Y_j^{r}f^q(X_j)\,|\,
{\cal F}_{j-1}]\,)^2\big]
\\
&\leq
2{\rm E}[Y_1^{2r}]{\textstyle\sum_{j=1}^{\infty}}\frac{1}{j^2}<+\infty.
\end{split}
\end{equation*}
Therefore, by Kronecher's lemma, we find that
\begin{equation*}
\frac{1}{n}{\textstyle\sum_{j=1}^n}
\big(Y_j^{r}f^q(X_j)-{\rm E}[Y_j^{r}f^q(X_j)\,|\, 
{\cal F}_{j-1}]\big)
\stackrel{a.s.}\longrightarrow 0.
\end{equation*}
On the other hand, since $Y_j$ is independent of 
${\cal F}_{j}^X\vee{\cal F}_{j-1}^Y$ by assumption,  
we have
\begin{equation*}
{\rm E}[Y_j^{r}f^q(X_j)\,|\, {\cal F}_{j-1}]=
{\rm E}[Y_1^r]{\rm E}[f^q(X_j)\,|\, {\cal F}_{j-1}]
={\rm E}[Y_1^r]V^{f^q}_{j-1}
\stackrel{a.s.}\longrightarrow 
{\rm E}[Y_1^r]V_{f^q}.
\end{equation*}
Since $n^{-1}\sum_{j=1}^{n}a_jd_j\stackrel{a.s.}\to ad$ whenever 
\begin{equation}\label{cond-serie}
a_j\geq 0,\quad d_j\stackrel{a.s.}\to d,\quad 
n^{-1}{\textstyle\sum_{j=1}^{n}}a_j\stackrel{a.s.}\to a,
\end{equation}
we obatin that  
\begin{equation*}
\frac{1}{n}{\textstyle\sum_{j=1}^n}
{\rm E}[Y_j^{r}f^q(X_j)\,|\, 
{\cal F}_{j-1}]
\stackrel{a.s.}\longrightarrow 
{\rm E}[Y_1^r]V_{f^q}
\end{equation*}
and so 
\begin{equation*}
\frac{1}{n}{\textstyle\sum_{j=1}^n} Y_j^{r}f^q(X_j)
\stackrel{a.s.}\longrightarrow 
{\rm E}[Y_1^r]V_{f^q}.
\end{equation*}
In particular, we get
\begin{equation*}
A^f_j\stackrel{a.s.}\longrightarrow \frac{V_f}{{\rm E}[Y_1]}.
\end{equation*}
Summing up, we have proved that $U^f_n$ is a sum of terms of the type
$n^{-1}\sum_{j=1}^{n}a_jd_j$, where $(a_j)$ and $(d_j)$ satisfy conditions 
(\ref{cond-serie}) and so we finally get that
$U^f_n$ converges a.s. to $U_f=\Delta(V_{f^2}-V_f^2)$.  By Theorem
\ref{clt2}, we conclude that $S_n^f$ converges $\cal A$-stably to
the Gaussian kernel ${\cal N}\big(0,\Delta(V_{f^2}-V_f^2)\big)$.}
$\qquad\Diamond$

\end{example}

The second result of this section is contained in the following
theorem.

\begin{theorem}\label{clt1}
{  Let $(X_n)_{n\geq 1}$ be a GOS and $f$ be a bounded Borel function. 
Using the previous notation, for $n\geq 0$ set
\begin{equation*}
Q_n:=p_{n+1,n+1}=1-\frac{r_{n+1}}{r_n}
\qquad\hbox{and}
\qquad
W^f_n:=\sqrt{n}(V^f_n-V_f).
\end{equation*}
Suppose that the following conditions are satisfied: 
\\[5pt] \indent
(i) $n\textstyle\sum_{k\geq n}
Q_k^2\stackrel{a.s.}\longrightarrow H$, where $H$ is a positive real
random variable. 
\\[5pt]\indent 
(ii) $\textstyle\sum_{k\geq 0} k^2\,{\rm E}[Q_k^4]<\infty$.
 \\[5pt] \indent Then the sequence $(W^f_n)_{n\geq 0}$
converges to the Gaussian kernel
 $${\cal N}\big(0,H(V_{f^2}-V_f^2)\,\big)$$ in the sense of the almost sure
 conditional convergence with respect to the filtrations 
$\CF=(\CF^X_n\vee\CF^Y_n)_{n\geq 0}$ and 
$\CG=(\CF^X_n\vee\CF^Y_{\infty})_{n\geq 0}$. 
\\[5pt]
In particular, we have
$$W^f_n\stackrel{{\cal A}-\hbox{stably}}\longrightarrow
{\cal N}\big(0,H(V_{f^2}-V_f^2)\,\big).
$$ 
}
\end{theorem}

\begin{cor}\label{cor-clt1}
  Using the notation of Theorem \ref{clt1}, let us set { for $k\geq 0$}
$$
\rho_k=\frac{1}{r_{k+1}}-\frac{1}{r_k}
$$
and assume the following conditions:
\\[5pt]
\indent (a)    
$r_k\leq c_k$ a.s. with  $\sum_{k\geq 0}k^2c_{k+1}^4<\infty$
and $kr_k\stackrel{a.s.}\to \alpha$,
where $c_k,\,\alpha$ are strictly positive constants.
\\[5pt]
\indent (b) The random variable $\rho_k$ are independent and
identically distributed with ${\rm E}[\rho_k^4]<\infty$.
\\[5pt]
Finally, let us set $\beta:={\rm E}[\rho_k^2]$ and $h:=\alpha^2\beta$.
\\[5pt]
\indent Then, 
 the conclusion of Theorem \ref{clt1} holds true with $H$ 
equal to the constant $h$.
\end{cor}

\begin{example}\label{example-clt1}
\rm Let us consider the setting of Example \ref{ex-2} with
$$
r_k=\frac{\th}{\th+\textstyle\sum_{i=1}^k Y_i}.
$$ where $\th>0$ and the random variables $Y_n$ are identically
distributed with $Y_n\geq\gamma>0$ and ${\rm E}[Y_n^4]<+\infty$.  Let us set
${\rm E}[Y_1]=m$ and ${\rm E}[Y_1^2]=\delta$. We have $r_{k}\leq
c_k=\th/(\th+\gamma k)$ and, by the strong law of large numbers, we have
\begin{equation*}
kr_{k}=\frac{\th k}{\th+\textstyle\sum_{i=1}^k Y_i}
\stackrel{a.s.}\longrightarrow \th/m.
\end{equation*}
Furthermore we have
\begin{equation*}
\rho_k=\frac{1}{r_{k+1}}-\frac{1}{r_k}=
\frac{Y_{k+1}}{\th}
\end{equation*}
and so $\beta={\rm E}[\rho_k^2]=\delta/\th^2$. Therefore
 the above corollary holds with $h=\delta/m^2$.$\qquad\Diamond$
\end{example}

The particular generalized P\`olya urn discussed in \cite{cri}
(Cor.~4.1) and in \cite{may-paganoni-secchi} is included in the above
example.

\section{Random partition induced by a GOS}\label{Sec.5}

Exchangeable species sampling
sequences are strictly connected with exchangeable random partitions.
Random partitions have been studied extensively, see, for
instance \cite{pitman} and the references theirin. 

In this section we investigate some properties of the length $L_n$ of
the random partition induced by a GOS at time $n$, i.e. the random
number of distinct values of GOS until time $n$.

Let $A_0:=E$ and
$A_n(\omega):=E\setminus\{X_1(\omega),\dots,X_n(\omega)\}= \{y\in
E:\;y\notin\{X_1(\omega),\dots,X_n(\omega)\}\}$ for $n\geq 1$ and 
define the following $\CF_n$-measurable random variable:
$$
s_n:=r_{n}(Y(n))\mu(A_n)=r_{n}\mu(A_n).
$$ 

\begin{remark}{\rm 
Reconsidering the species interpretation, given $X(n)=(X_1, \dots X_n)$ 
and $Y(n)=(Y_1,\dots,Y_n)$, the species
of the $(n+1)$-th individual is a new species with probability $s_n$
and one of the species observed so far with probability $1-s_n$.
In particular one has
\begin{equation*}
P[L_{n+1}=L_n+1\,|\,{\cal F}_{n}]=s_n=r_{n}\mu(A_n). 
\end{equation*} 
If the probability measure $\mu$ is diffuse, then $s_n=r_n$.  }
\end{remark}

If $\mu$ is diffuse and the coefficients $r_n$ are deterministic
( such as in Example \ref{ex-2bis}), then
the sequence of the increments $(L_n-L_{n-1})_{n\geq 1}$ (with $L_0:=0$)
is a sequence of independent random variables such that, for each $n$,
the distribution of $L_n-L_{n-1}$ is a Bernoulli distribution with
parameter $r_{n-1}$, hence it is immediate to deduce, under suitable
conditions, both a strong law of large numbers and a central limit
theorem  for $(L_n)_{n \geq 1}$.  

In this section we prove a law of large numbers and a central 
limit theorem for a GOS. Moreover, some examples of GOS that satisfy the 
hypotheses of these results are given.

\begin{theorem} \label{slln} Let $(X_n)_{n\geq 1}$ be a generalized
  species sampling sequence.  Suppose that there exists a sequence
  $(h_n)_{n\geq 1}$ of real numbers and a random variable $L$ such
  that the following properties hold:
\begin{equation*}
h_n\geq 0,\quad h_n\uparrow +\infty,\quad
{\textstyle\sum_{j\geq 1}}\frac{{\rm E}[s_j(1-s_j)]}{h_j^2}<+\infty,
\quad 
\frac{1}{h_n}{\textstyle\sum_{j=0}^{n}}s_{j}\stackrel{a.s.}\longrightarrow L.
\end{equation*}
Then we have ${L_n}/{h_n}\stackrel{a.s.}\longrightarrow L.$
\end{theorem}

\begin{remark}{\rm
\rm  Let us note that, for each $n$, we have
\begin{equation*} 
{\rm E}[L_{n+1}\,|\,{\cal F}_n]=L_n+s_n\geq L_n.  
\end{equation*} 
Hence the sequence $(L_n)_{n\geq 0}$ is
 a positive submartingale with ${\rm E}[L_{n+1}]={\rm
 E}[L_n]+{\rm E}[s_n]$. Therefore $(L_n)_{\geq 0}$ is bounded
in $L^1$ if and only if we have $\sum_{k\geq 0} {\rm E}[s_k]<
+\infty$ and, in this case, $(L_n)_{n\geq 0}$ converges almost surely to
an integrable random variable. It follows that, for each
sequence $(h_n)_{n\geq 0}$ with $h_n\to +\infty$, the ratio $L_n/h_n$
goes almost surely to zero. An example of this situation is given by 
Example \ref{ex-2tris} with $\mu$ diffuse. 
Indeed, in this case, we have ${\rm E}[s_n]=
{\rm E}[Y_n]=(1/2)^{n-1}$.
} 
\end{remark}

\begin{theorem}\label{clt}
  Let $(X_n)_{n\geq 1}$ be a GOS with $\mu$ diffuse and suppose there
  exists a sequence $(h_n)_{n\geq 1}$ of real numbers and a positive
  random variable $\sigma^2$ such that the following properties hold:
\begin{equation*}
h_n\geq 0,\quad h_n\uparrow +\infty,\quad
\sigma_n^2:=\frac{\textstyle\sum_{j=1}^n r_j(1-r_j)}{h_n}
\stackrel{a.s.}\longrightarrow\sigma^2.
\end{equation*}
Then, setting  $R_n:=\sum_{j=1}^n r_j$, we have 
\begin{equation*}
T_n:=\frac{L_n-R_{n-1}}{\sqrt{h_n}}
\stackrel{{\cal A}-\hbox{stably}}\longrightarrow{\cal N}(0,\sigma^2).
\end{equation*}
\end{theorem}

\begin{cor}
Under the same assumptions of Theorem \ref{clt}, if
$P(\sigma^2>0)=1$, then we have
\begin{equation*}
\frac{T_n}{\sigma_n}=\frac{(L_n-R_{n-1})}
{\sqrt{\textstyle\sum_{j=1}^n r_j(1-r_j)}}
\stackrel{{\cal A}-\hbox{stably}}\longrightarrow{\cal N}(0,1).
\end{equation*}
\end{cor}

\begin{example} \label{cinese-potenza}
  \rm 
 Let us consider Example \ref{ex-2bis} with $\mu$ diffuse and  
$$a_n=\frac{\th}{\th+n^{1-\alpha}}$$
with $\th>0$ and $0<\alpha< 1$. We have $s_n=r_n=a_n$
and, setting
$h_n=n^{\alpha}$ and $L=\th/\alpha$, the assumptions of Theorem
\ref{slln} are satisfied. Indeed we have
\begin{equation*}
\begin{split}
{\textstyle\sum_{j}}\frac{r_j(1-r_j)}{j^{2\alpha}}
&=
{\textstyle\sum_{j}}\frac{j^{1-\alpha}}{(\th+j^{1-\alpha})^2j^{2\alpha}}
\\
&=
{\textstyle\sum_{j}}
\left(\frac{j^{1-\alpha}}{\th+j^{1-\alpha}}\right)^2\frac{1}{j^{\alpha+1}}
<+\infty.
\end{split}
\end{equation*}
Moreover, since
\begin{equation}\label{power}
\frac{1}{n^{\alpha}}{\textstyle\sum_{j=1}^{n}} \frac{1}{j^{1-\alpha}}
\longrightarrow\frac{1}{\alpha}\qquad\hbox{for } \alpha\in(0,1),
\end{equation}
we have
$$
\frac{1}{n^{\alpha}}R_n=
\frac{1}{n^{\alpha}}{\textstyle\sum_{j=1}^{n}} \frac{\th}{\th+j^{1-\alpha}}
\longrightarrow
\frac{\th}{\alpha}.
$$
Thus we have $L_n/n^{\alpha}\stackrel{a.s.}\longrightarrow \th/\alpha$.
Finally, since
\begin{equation}\label{res1}
\frac{1}{h_n}{\textstyle\sum_{j=1}^n} a_jb_j\to b,
\end{equation}
provided that $a_j\geq 0$,  $\sum_{j=1}^n a_j/h_n \to 1$
and $b_n \to b$ as $n \to +\infty$, 
it is easy to see that
\begin{equation*}
\sigma_n^2=\frac{\textstyle\sum_{j=1}^n r_j(1-r_j)}{n^{\alpha}}
=\frac{\th}{n^{\alpha}}
{\textstyle\sum_{j=1}^n}\frac{j^{1-\alpha}}{(\th+j^{1-\alpha})^2}
=\frac{\th}{n^{\alpha}}
{\textstyle\sum_{j=1}^n}\left(\frac{j^{1-\alpha}}{\th+j^{1-\alpha}}\right)^2
\frac{1}{j^{1-\alpha}}
\to \th/\alpha.
\end{equation*}
Therefore, by Theorem \ref{clt}, we obtain 
\begin{equation*}
T_n=\frac{L_n-R_{n-1}}{n^{\alpha/2}}
\stackrel{\cal A-\hbox{stably}}\longrightarrow{\cal N}(0,\th).
\end{equation*}
$\qquad\Diamond$

\end{example}

\begin{example} \label{cinese-random} 
\rm Let us consider the setting of Example \ref{ex-2} with $\mu$
diffuse and
$$
r_n=\frac{\th}{\th+\textstyle\sum_{i=1}^n Y_i}.
$$ where $\th>0$ and the random variables $Y_n$ are independent identically
distributed positive random variable with ${\rm E}[Y_n]=m>0$. Then
$s_n=r_n$ and, setting $h_n=\log n$ and $L=c/m$, the assumptions
of Theorem \ref{slln} are satisfied. Indeed 
$$
{\textstyle\sum_{j}}\frac{{\rm E}\big[r_j(1-r_j)\big]}{(\log j)^2}
\leq
{\textstyle\sum_{j}}\frac{1}{(\log j)^2}<+ \infty.
$$
Moreover, by the strong law of large numbers, we have
$$
\left(\frac{\th}{j}+
\frac{1}{j}{\textstyle\sum_{i=1}^j}Y_i\right)^{-1}
\stackrel{a.s.}\longrightarrow 1/m.
 $$ 
Therefore, 
since $\frac{1}{\log n}\sum_{j=1}^n \frac{1}{j}\to 1$, 
by (\ref{res1}),  we can conclude that 
$$
\frac{1}{\log n}R_n
=
\frac{\th}{\log n}
{\textstyle\sum_{j=1}^n}\frac{1}{\th+{\textstyle\sum_{i=1}^j}Y_i}
= 
\frac{\th}{\log n}
{\textstyle\sum_{j=1}^n} 
\frac{1}{j}\left(\frac{\th}{j}+
\frac{1}{j}{\textstyle\sum_{i=1}^j}Y_i\right)^{-1}
\stackrel{a.s.}\longrightarrow \frac{\th}{m}
$$
and so $L_n/\log n\stackrel{a.s.}\longrightarrow \th/m$.  
Moreover,
by (\ref{res1}) and the strong law of large numbers, we have
\begin{equation*}
\begin{split}
\sigma_n^2&=\frac{\textstyle\sum_{j=1}^n r_j(1-r_j)}{\log n}
=\frac{\th}{\log n}{\textstyle\sum_{j=1}^n} 
\frac{\textstyle\sum_{i=1}^j Y_i}{(\theta+\textstyle\sum_{i=1}^j Y_i)^2}\\
&=\frac{\th}{\log n}
{\textstyle\sum_{j=1}^n }\left(\frac{\textstyle\sum_{i=1}^j Y_i/j}
{\th/j+\textstyle\sum_{i=1}^j Y_i/j}\right)^2
\frac{j}{\textstyle\sum_{i=1}^j Y_i}
\,\frac{1}{j}
\to \th/m.
\end{split}
\end{equation*}
Therefore, by Theorem \ref{clt}, we obtain 
\begin{equation*}
T_n=\frac{L_n-R_{n-1}}{\sqrt{\log n}}
\stackrel{\cal A-\hbox{stably}}\longrightarrow{\cal N}(0,\th/m)
\end{equation*}
and so
\begin{equation*}
\frac{L_n-R_{n-1}}{\sqrt{\frac{\th}{m}\log n}}
\stackrel{\cal A-\hbox{stably}}
\longrightarrow {\cal N}(0,1).
\end{equation*}
If we take $Y_i=1$ for all $i$, we find the well known results for the
asymptotic distribution of the length of the random partition obtained
with the Blackwell--McQueen urn scheme. Indeed, since
\( \sum_{j=1}^{n} {j}^{-1} -\log n=\gamma+O(\frac{1}{n})\),
one gets
\begin{equation*}
\frac{L_n-\th\log n}{\sqrt{\th\log n}}\stackrel{\cal A-\hbox{stably}}
\longrightarrow {\cal N}(0,1).
\end{equation*}
See,  for instance, pages 68-69 in \cite{pitman}.$\qquad\Diamond$
\end{example}

\section{Proofs.}

This section contains all the proofs of the paper. Recall that
$$
\CF_n=\CF^X_n\vee\CF^Y_n\qquad\hbox{and}\qquad
\CG_n=\CF^X_n\vee\CF^Y_{\infty}=\CF_n\vee\sigma(Y_{n+j}:j\geq 1)
$$
and so condition $(h3)$ of the definition of generalized species
sampling sequence implies that 
$$
V_n^g:={\rm E}[g(X_{n+1})\,|\,{\CG}_n]={\rm E}[g(X_{n+1})\,|\,{\cal F}_n]
$$
for each bounded Borel function $g$ on $E$.

\subsection{Proof of  Theorem \ref{cond-cidness}}

We start with a useful lemma.

\begin{lemma}\label{lemma-prob-cond} 
  If $(X_n)_{n\geq 1}$ is a generalized species sampling sequence,
then we have
\begin{equation*}
P[n+1\in\pi_l^{(n+1)}\,|\,{\CG}_n]=
P[X_{n+1}=X_l^*\,|\,{\CF}_n]=
{\textstyle\sum_{j\in\pi_l^{(n)}}}p_{n,j}(\pi^{(n)},Y(n))
+r_n(\pi^{(n)},Y(n))\mu(\{X^*_l\})
\end{equation*}
for each $l=1,\dots,L_n$. 
Moreover, 
\begin{equation*}
{\rm E}[I_{\{L_{n+1}=L_n+1\}}f(X_{n+1})\,|\,\CG_{n}]=
{\rm E}[I_{\{L_{n+1}=L_n+1\}}f(X_{n+1})\,|\,\CF_{n}]=
r_n(\pi^{(n)},Y(n))\textstyle\int_{A_n}f(y)\,\mu({\rm d}y).
\end{equation*} 
holds true with $A_0:=E$ and $A_n$ the random ``set'' defined by
$$
A_n(\omega):=E\setminus\{X_1(\omega),\dots,X_n(\omega)\}= \{y\in
E:\;y\notin\{X_1(\omega),\dots,X_n(\omega)\}\}
\quad\hbox{for } n\geq 1.
$$
In particular, we have 
\begin{equation*}
P[L_{n+1}=L_n+1\,|\,\CG_{n}]=P[L_{n+1}=L_n+1\,|\,{\cal F}_{n}]=
r_{n}(\pi^{(n)},Y(n))\mu(A_n):=s_n(\pi^{(n)},Y(n)) 
\end{equation*} 
If $\mu$ is diffuse, we have
\begin{equation*}
P[n+1\in\pi_l^{(n+1)}\,|\,{\CG}_n]=
{\rm P}[X_{n+1}=X_l^*\,|\,{\CF}_n]=
{\textstyle\sum_{j\in\pi_l^{(n)}}}p_{n,j}(\pi^{(n)},Y(n))
\end{equation*}
for each $l=1,\dots,L_n$ and 
\begin{equation*}
{\rm E}[I_{\{L_{n+1}=L_n+1\}}f(X_{n+1})\,|\,\CG_{n}]=
{\rm E}[I_{\{L_{n+1}=L_n+1\}}f(X_{n+1})\,|\,\CF_{n}]=
r_n(\pi^{(n)},Y(n)){\rm E}[f(X_1)]
\end{equation*} 
and 
\begin{equation*}
P[L_{n+1}=L_n+1\,|\,\CG_{n}]=P[L_{n+1}=L_n+1\,|\,{\cal F}_{n}]
=r_{n}(\pi^{(n)},Y(n)). 
\end{equation*} 
\end{lemma}

{\bf Proof.} 
Since $\CG_n=\CF_n\vee\sigma(Y_{n+j}:j\geq 1)$, condition $(h_3)$ implies 
that
\begin{equation*}
P[n+1\in\pi_l^{(n+1)}\,|\,{\CG}_n]=
P[X_{n+1}=X_l^*\,|\,{\CG}_n]
P[X_{n+1}=X_l^*\,|\,{\CF}_n].
\end{equation*}
Hence, by assumption $(h2)$, we have
\begin{equation*}
\begin{split}
P[X_{n+1}=X_l^*\,|\,{\CF}_n]&=
\textstyle\sum_{i=1}^n p_{n,i}(\pi^{(n)},Y(n))\delta_{X_i}(X^*_l)
+r_n(\pi^{(n)},Y(n))\mu(\{X^*_l\})\\
&=
{\textstyle\sum_{j\in\pi_l^{(n)}}} p_{n,j}(\pi^{(n)},Y(n))
+r_n(\pi^{(n)},Y(n))\mu(\{X^*_l\}).
\end{split}
\end{equation*}
for each $l=1,\dots,L_n$. 
If $\mu$ is diffuse, we obtain
\begin{equation*}
P[X_{n+1}=X_l^*\,|\,{\CF}_n]=
{\textstyle\sum_{j\in\pi_l^{(n)}}} p_{n,j}(\pi^{(n)},Y(n))
\end{equation*}
for each $l=1,\dots,L_n$. 

Now, we observe that
$$
I_{\{L_{n+1}=L_n+1\}}=I_{B_n}(X_1,\dots,X_n,X_{n+1})
$$ where
$B_n=\{(x_1,\dots,x_{n+1}):\;x_{n+1}\notin\{x_1,\dots,x_n\}\}$. Thus,
by $(h_3)$ and $(h_2)$, we have 
\begin{equation*}
\begin{split}
{\rm E}[I_{\{L_{n+1}=L_n+1\}}f(X_{n+1})\,|\,\CG_{n}]&=
{\rm E}[I_{\{L_{n+1}=L_n+1\}}f(X_{n+1})\,|\,\CF_{n}]\\
&=
\textstyle\int I_{B_n}(X_1,\dots,X_n,y)f(y)K_{n+1}(\cdot,{\rm d}y)\\
&=
\textstyle\sum_{i=1}^n p_{n,i}(\pi^{(n)},Y(n))
\int_{A_n}f(y)\delta_{X_i}({\rm d}y)
+r_n(\pi^{(n)},Y(n))\int_{A_n}f(y)\mu({\rm d}y)\\
&=r_n(\pi^{(n)},Y(n))\textstyle\int_{A_n}f(y)\mu({\rm d}y).
\end{split}
\end{equation*}
If we take $f=1$, we get
\begin{equation*}
P[U_{n+1}=1\,|\,\CG_{n}]=P[U_{n+1}=1\,|\,{\cal F}_{n}]=
r_n(\pi^{(n)},Y(n))\mu(A_n).
\end{equation*}
Finally, if $\mu$ is diffuse, then $\mu(A_n(\omega))=1$ for each
$\omega$ and so we have
\begin{equation*}
\textstyle\int_{A_n}f(y)\mu({\rm d}y)=
{\rm E}[f(X_1)].
\end{equation*}
\qed

{\bf Proof of Theorem \ref{cond-cidness}. }   
   Let us fix a bounded Borel function $f$
on $E$.  Using the given prediction rule, we have 
\begin{equation*}
\begin{split}
V^f_n&=
{\textstyle\sum_{i=1}^n} 
p_{n,i}(\pi^{(n)},Y(n))f(X_i)+r_n(\pi^{(n)},Y(n)){\rm E}[f(X_1)]\\
&=
{\textstyle\sum_{j=1}^{L_n}} p^*_{n,j}(\pi^{(n)})f(X_j^*)+r_n{\rm E}[f(X_1)].\\
\end{split}
\end{equation*}
The sequence $(X_n)$ is $\CG$-cid if and only if for each   bounded
Borel function $f$ on $E$, the sequence $(V^f_n)_{n\geq 0}$ is a
$\CG$-martingale. We observe that we have (for the sake of
simplicity we skip the dependence on $(Y_n)_{n\geq 1}$)
\begin{equation*}
\begin{split}
{\rm E}[V^f_{n+1}\,|\,{\CG}_n]&=
{\textstyle\sum_{i=1}^n} 
f(X_i)E_i+{\rm E}[p_{n+1,n+1}(\pi^{(n+1)})f(X_{n+1})\,|\,{\CG}_n]+
{\rm E}[r_{n+1}\,|\,{\CG}_n]\bar f\\
&=
{\textstyle\sum_{j=1}^{L_n} f(X_j^*)\sum_{i\in\pi_j^{(n)}}}E_i+
{\rm E}[p_{n+1,n+1}(\pi^{(n+1)})f(X_{n+1})\,|\,{\CG}_n]+
{\rm E}[r_{n+1}\,|\,{\CG}_n]\bar f
\end{split}
\end{equation*}
where $E_i={\rm E}[p_{n+1,i}(\pi^{(n+1)})\,|\,{\CG}_n]$ and
$\bar f={\rm E}[f(X_1)]$.\\
Now we are going to compute the various conditional expectations which
appear in the second member of above equality.  Since $\mu$ is
  diffuse, using Lemma \ref{lemma-prob-cond}, we have 

\begin{equation*}
\begin{split}
E_i&={\rm E}[p_{n+1,i}(\pi^{(n+1)})\,|\,{\CG}_n]\\
&=
\textstyle{\sum_{l=1}^{L_n}}
{\rm E}[ I_{\{n+1\in\pi_l^{(n+1)}\}} p_{n+1,i}(\pi^{(n+1)})\,|\,{\CG}_n]+
{\rm E}[I_{\{L_{n+1}={L_n}+1\}}p_{n+1,i}(\pi^{(n+1)})\,|\,{\CG}_n]\\
&=
\textstyle{\sum_{l=1}^{L_n}}
p_{n+1,i}([\pi^{(n)}]_{l^+}){\rm E}[ I_{\{n+1\in\pi_l^{(n+1)}\}}\,|\,{\CG}_n]
+{\rm E}[I_{\{L_{n+1}={L_n}+1\}}\,|\,{\CG}_n] p_{n+1,i}([\pi^{(n)};n+1])\\
&=
{\textstyle\sum_{l=1}^{L_n} }
p_{n+1,i}([\pi^{(n)}]_{l^+})
{\textstyle\sum_{j\in\pi_l^{(n)}}}p_{n,j}(\pi^{(n)})
+r_n p_{n+1,i}([\pi^{(n)};n+1])\\
&=
{\textstyle\sum_{l=1}^{L_n}} 
p_{n+1,i}([\pi^{(n)}]_{l^+}) p_{n,l}^*(\pi^{(n)})+
r_n p_{n+1,i}([\pi^{(n)};n+1])
\end{split}
\end{equation*}
and so
\begin{equation*}
\begin{split}
\sum_{i\in\pi_j^{(n)}}E_i&=
\sum_{l=1,l\neq j}^{L_n} p^*_{n+1,j}([\pi^{(n)}]_{l^+})p^*_{n,l}(\pi^{(n)})
+\sum_{i\in\pi_j^{(n)}} p_{n+1,i}([\pi^{(n)}]_{j^+})p^*_{n,j}(\pi^{(n)})
+r_n p^*_{n+1,j}([\pi^{(n)};n+1])\\
&=
{\textstyle \sum_{l=1}^{L_n}} 
p^*_{n+1,j}([\pi^{(n)}]_{l^+})p^*_{n,l}(\pi^{(n)})
-p_{n+1,n+1}([\pi^{(n)}]_{j^+} )p^*_{n+1,j}(\pi^{(n)})
+r_np^*_{n+1,j}([\pi^{(n)};n+1])
\end{split}
\end{equation*}
Moreover, using Lemma \ref{lemma-prob-cond} again, we have
\begin{equation*}
\begin{split}
&{\rm E}[p_{n+1,n+1}(\pi^{(n+1)})f(X_{n+1})\,|\,{\CG}_n]=\\
&\sum_{l=1}^{L_n}
{\rm E}[I_{\{n+1\in\pi_l^{(n+1)}\}}p_{n+1,n+1}(\pi^{(n+1)})f(X_{n+1})
\,|\,{\CG}_n]+
{\rm E}[I_{\{L_{n+1}={L_n}+1\}}p_{n+1,n+1}(\pi^{(n+1)})f(X_{n+1})
\,|\,{\CG}_n]=\\
&
\sum_{l=1}^{L_n}
{\rm E}[I_{\{n+1\in\pi_l^{(n+1)}\}}\,|\,{\CG}_n]
p_{n+1,n+1}([\pi^{(n)}]_{l+})f(X^*_{l})+
{\rm E}[I_{\{L_{n+1}={L_n}+1\}}f(X_{n+1})\,|\,{\CG}_n]
p_{n+1,n+1}([\pi^{(n)}];n+1)=\\
&
{\textstyle\sum_{l=1}^{L_n}}
\left({\textstyle\sum_{k\in\pi_l^{(n)}}}p_{n,k}(\pi^{(n)})\right) 
p_{n+1,n+1}([\pi^{(n)}]_{l^+})f(X_l^*)
+r_n p_{n+1,n+1}([\pi^{(n)}];n+1)\bar f=\\
&
{\textstyle\sum_{l=1}^{L_n}}
p^*_{n,l}(\pi^{(n)}) p_{n+1,n+1}([\pi^{(n)}]_{l^+})f(X_l^*)
+r_n p_{n+1,n+1}([\pi^{(n)}];n+1)\bar f.\\
\end{split}
\end{equation*}
Finally we have
\begin{equation*}
\begin{split}
{\rm E}[r_{n+1} \,|\,{\CG}_n]&= 
1 - \textstyle\sum_{i=1}^{n+1}{\rm E}[p_{n+1,i}(\pi^{(n+1)})\,|\,{\CG}_n]\\
&=1 - 
\textstyle\sum_{i=1}^n E_i - E_{n+1}\\
&=1 - 
{\textstyle\sum_{i=1}^n} E_i-
{\textstyle \sum_{l=1}^{L_n}}
p^*_{n,l}(\pi^{(n)})p_{n+1,n+1}([\pi^{(n)}]_{l^+})
-r_n p_{n+1,n+1}([\pi^{(n)}];n+1)
\end{split}
\end{equation*}
Thus we get 
\begin{equation*}
{\rm E}[V^f_{n+1}\,|\,{\CG}_n]=
\textstyle\sum_{j=1}^{L_n} c_{n,j} f(X_j^*) + 
(1-\sum_{j=1}^{L_n} c_{n,j})\bar f
\end{equation*}
where
\begin{equation*}
\begin{split}
c_{n,j}&=
{\textstyle\sum_{i\in\pi_j^{(n)}}}E_i+
p_{n+1,n+1}([\pi^{(n)}]_{j^+})p^*_{n,j}(\pi^{(n)})\\
&=r_n p^*_{n+1,j}([\pi^{(n)};n+1])+
{\textstyle\sum_{l=1}^{L_n}}
p^*_{n+1,j}([\pi^{(n)}]_{l^+})p^*_{n,l}(\pi^{(n)})\\
\end{split}
\end{equation*}
{We can conclude that $(X_n)_{n\geq 1}$ is $\CG$-cid if and only if we have, 
for each bounded Borel function $f$ on $E$ and each $n$
\begin{equation*}
\textstyle\sum_{j=1}^{L_n} p^*_{n,j}f(X_j^*)+r_n\bar f
=\textstyle\sum_{j=1}^{L_n} c_{n,j} f(X_j^*) + 
(1-\sum_{j=1}^{L_n} c_{n,j}) \bar f
\qquad P\hbox{-almost surely.}
\end{equation*}
}
Since $E$ is a Polish space, we may affirm that $(X_n)_{n\geq 1}$ is $\CG$-cid
if and only if, for each $n$, we have $P$-almost surely
\begin{equation*}
\textstyle\sum_{j=1}^{L_n} p^*_{n,j}\delta_{X_k^*}(\cdot)+r_n\mu(\cdot)
=\textstyle\sum_{j=1}^{L_n} c_{n,j} \delta_{X_k^*}(\cdot) + 
(1-\sum_{j=1}^{L_n} c_{n,j}) \mu(\cdot)
\end{equation*}
But this last equality holds if and only if, for each $n$, we have
$P$-almost surely
\begin{equation*}
p^*_{n,j}=c_{n,j}
\qquad\hbox{for } 1\leq j\leq L_n\; ;
\end{equation*}
that is 
\begin{equation*}
p_{n,j}^*(\pi^{(n)})=
r_n p_{n+1,j}^*([\pi^{(n)};\{n+1\}])
+
{\textstyle\sum_{l=1}^{L_n}} 
p_{n+1,j}^*([\pi^{(n)}]_{l+},)p_{n,l}^*(\pi^{(n)})
\end{equation*}
This is exactly the condition in the statement of the Theorem
\ref{cond-cidness}.

\qed

\subsection{Proofs of Section \ref{Sec.4}}

{\bf Proof of Theorem \ref{clt2}. }
We will use Theorem \ref{fam-tri} in the Appendix. For each $n\geq 1$,
let us set
\begin{equation*}
D^f_{n}=\sqrt{n} (M^f_{n}- V_f)=
\frac{1}{\sqrt{n}}
\big[{\textstyle\sum_{k=1}^n} f(X_k)-nV_f\big],
\end{equation*}
and, for $0\leq j\leq n$,
\begin{equation*}
L^f_{n,j}={\rm E}[D^f_{n}\,|\,\CG_j]
\qquad
{\cal F}_{n,j}=\CG_j.
\end{equation*}
Then, for each $n\geq 1$, the sequence $(L_{n,j})_{0\leq j\leq n}$
is a martingale with respect to $({\cal F}_{n,j})_{0\leq j\leq n}$ such that
$L_{n,0}={\rm E}[D^f_n|\CG_0]=0$ and, for $1\leq j\leq n$,
\begin{equation*}
L^f_{n,j}-L^f_{n,j-1}=
{\rm E}[D^f_{n}\,|\,\CG_j]-{\rm E}[D^f_{n}\,|\,\CG_{j-1}]
=Z^f_{n,j}.
\end{equation*}
Indeed, using (\ref{legame-fondamentale}) we have
\begin{equation*}
\begin{split}
&{\rm E}[D^f_{n}\,|\,\CG_j]-{\rm E}[D^f_{n}\,|\,\CG_{j-1}]
\\
&=
\frac{1}{\sqrt{n}}
\big[{\textstyle
\sum_{k=1}^{j}} f(X_k) +(n-j)V^f_j-nV^f_j-
\textstyle{\sum_{k=1}^{j-1}} f(X_k) -(n-j+1)V^f_{j-1}+nV^f_{j-1}
\big]\\
&=\frac{1}{\sqrt{n}}
\big[f(X_j) -jV^f_j+(j-1)V^f_{j-1}\big].
\end{split}
\end{equation*}
\\
Moreover, we have 
$$
S^f_n={\rm E}[D^f_{n}\,|\,\CG_n]=L^f_{n,n}=
\textstyle\sum_{j=1}^n Z^f_{n,j}.$$ 
Finally, we have 
\begin{equation*}
{\cal H}_j=\textstyle\liminf_n{\cal F}_{n,j\wedge n}=
\textstyle\liminf_n\CG_{j\wedge n}=\CG_{j}
\end{equation*}
and, if we set
\begin{equation*}
{\cal H}={\textstyle\bigvee_{j\geq 0}}{\cal H}_j
={\textstyle\bigvee_{j\geq 0}}\CG_j,
\end{equation*}
then the random variable $U_f$ is measurable with respect to the
$\sigma$-field ${\cal H}\vee{\cal N}$. 
At this point we can apply
Theorem \ref{fam-tri} and the proof of the first assertion is
concluded.
\\[5pt]
\indent If conditions (a1) holds, then condition (a) is obviously
verified. Moreover we have
\begin{equation*}
Z^f_{n,j}=\frac{1}{\sqrt{n}}\, Z^f_j
\end{equation*}
where 
$$
Z^f_j=f(X_j)-jV^f_j+(j-1)V^f_{j-1}.
$$
We can write
\begin{equation*}
\frac{1}{n}(Z^f_{n})^2
=(Z^f_{n,n})^2=
{\textstyle\sum_{j=1}^{n}} (Z^f_{n,j})^2-
\frac{1}{n}{\textstyle\sum_{j=1}^{n-1}}(Z^f_{j})^2=
U^{f}_{n}-\frac{n-1}{n}U^{f}_{n-1}
\stackrel{a.s.}\longrightarrow 0,
\end{equation*}
This fact implies that 
$$(Z^f_n)^*={\textstyle
\sup_{1\leq j\leq n}}|Z^f_{n,j}|
\stackrel{a.s.}\longrightarrow 0,$$
Indeed,
\begin{equation*}
{\textstyle\sup_{0\leq j\leq n}} (Z^f_{n,j})^2=
\frac{1}{n}\,{\textstyle\sup_{0\leq j\leq n}} (Z^f_j)^2
\stackrel{a.s.}\longrightarrow 0.
\end{equation*}
Further, we have
\begin{equation*}
\begin{split}
{\rm E}\big[\big((X^f_n)^*\big)^2\big]&=
{\rm E}[\textstyle\sup_{1\leq j\leq n}(Z^f_{n,j})^2]
\leq\textstyle\sum_{j=1}^n{\rm E}[(Z^f_{n,j})^2]
=\textstyle\sum_{j=1}^n
{\rm E}\big[
\big(L^f_{n,j}-L^f_{n,j-1}\big)^2 
\big]\\
&=\textstyle\sum_{j=1}^n
{\rm E}\big[(L^f_{n,j})^2\big]-{\rm E}\big[(L^f_{n,j-1})^2\big]
={\rm E}\big[(L^f_{n,n})^2\big]={\rm E}[(S^f_{n})^2].
\end{split}
\end{equation*}
>From (b1) and the above relations, we obtain that the sequence
$\big((Z^f_n)^*\big)_n$ is bounded in $L^2$ and so we get condition (b).  
\qed

{\bf Proof of Theorem \ref{clt1}.}  Without loss of generality, we may
assume $|f|\leq 1$. It will be sufficient to prove that the sequence
$(V^f_n)_{n\geq 0}$ satisfies conditions (a) and (b) of Theorem
\ref{th-appendice}, with $U=H(V_{f^2}-V_f^2)$.  To this end, we
observe firstly that, after some calculations, we have
\begin{equation}\label{eq1}
V^f_k-V^f_{k+1}=\big[V^f_k-f(X_{k+1})\big]\,Q_{k}.
\end{equation}
>From this equality we get $|V^f_k-V^f_{k+1}|\leq Q_k$, and so, using
assumption (ii), we find
\begin{equation*}
\textstyle\sup_k k^2\,|V^f_k-V^f_{k+1}|^4\leq 
\textstyle\sum_{k\geq 0} k^2 Q_k^4\in L^1.
\end{equation*}
Furthermore, by (\ref{eq1}), we have 
\begin{equation*}
\textstyle\sum_{k\geq n}(V^f_k-V^f_{k+1})^2
=
\textstyle\sum_{k\geq n}\big[V^f_k-f(X_{k+1})\big]^2\,Q_k^2
\qquad\hbox{for }\; n\to +\infty.
\end{equation*}
Therefore, in order to complete the proof, it suffices to prove, for
$n\to+\infty$, the following convergence:
\begin{equation*}
n\textstyle\sum_{k\geq n}\big[V^f_k-f(X_{k+1})\big]^2\, Q_k^2
\stackrel{a.s.}\longrightarrow H(V_{f^2}-V_f^2).
\end{equation*}
The above convergence can be rewritten as
\begin{equation}\label{eq2}
n\textstyle\sum_{k\geq n}\big[(V^f_k)^2+f^2(X_{k+1})-2V^f_kf(X_{k+1})\big]
\,Q_k^2
\stackrel{a.s.}\longrightarrow H(V_{f^2}-V_f^2).
\end{equation}
Now, by assumption (i) and the almost sure convergence of $(V^f_k)_k$ to
$V_f$ and of $(V^{f^2}_k)_k$ to $V_{f^2}$, we have
\begin{equation}\label{eq3prima}
n\textstyle\sum_{k\geq n}V^f_k\, Q_k^{2}
\stackrel{a.s.}\longrightarrow V_fH
\end{equation}
\begin{equation}\label{eq3}
n\textstyle\sum_{k\geq n}(V^f_k)^2\, Q_k^{2}
\stackrel{a.s.}\longrightarrow (V_f)^2H
\end{equation}
\begin{equation}\label{eq4}
n\textstyle\sum_{k\geq n}V^{f^2}_k\, Q_k^{2}
\stackrel{a.s.}\longrightarrow V_{f^2}H
\end{equation}
Thus, it will be enough to prove the following convergence:
\begin{equation}\label{eq6}
n\textstyle\sum_{k\geq n}\big[g(X_{k+1})-V^g_k\big]\,Q_k^2
\stackrel{a.s.}\longrightarrow 0
\end{equation}
where $g$ is a bounded Borel function with $|g|\leq 1$. 
Indeed, from (\ref{eq6})  with $g=f^2$ and (\ref{eq4}), we obtain
\begin{equation}\label{eq7}
n\textstyle\sum_{k\geq n}f^2(X_{k+1})\, Q_k^2
\stackrel{a.s.}\longrightarrow V_{f^2}H,
\end{equation}
Moreover, from (\ref{eq6})  with $g=f$ and (\ref{eq3prima}), we obtain
\begin{equation}\label{eq8}
n\textstyle\sum_{k\geq n}f(X_{k+1})\, Q_k^2
\stackrel{a.s.}\longrightarrow V_{f}H, 
\end{equation}
and so, by the almost sure convergence of $(V^f_k)_k$ to
$V_f$, we get
\begin{equation}\label{eq9}
n\textstyle\sum_{k\geq n}V^f_k\, f(X_{k+1})\, Q_k^2
\stackrel{a.s.}\longrightarrow (V_{f})^2H, 
\end{equation} 
Then convergence relations (\ref{eq3}), (\ref{eq7}) and (\ref{eq9})
lead us to the desired relation (\ref{eq2}).  
\\[5pt] 
\indent In
order to prove (\ref{eq6}), we consider the process $(Z_n)_{n\geq 0}$
defined by
\begin{equation*}
Z_n:=\textstyle\sum_{k=0}^{n-1}k\,\big[g(X_{k+1})-V^g_k\big]\,Q_k^2.
\end{equation*}
It is a martingale with respect to the filtration $\CG=({\CG}_n)_{n\geq 0}$. 
Moreover, by assumption (ii), we have
\begin{equation}\label{limitata}
{\rm E}[Z_n^2]
=
\textstyle\sum_{k=0}^{n-1}k^{2}\,{\rm E}\big[(g(X_{k+1})-V^g_k)^2\,Q_k^4\big]
\leq
\textstyle\sum_{k\geq 0}k^{2}\,{\rm E}[Q_{k}^4]
<\infty.
\end{equation}
The martingale $(Z_n)_{n\geq 1}$ is thus bounded in $L^2$ and so it converges
almost surely; that is, the series
\begin{equation*}
\textstyle\sum_{k\geq 0}k\,\big[g(X_{k+1})-V^g_k\big]\,Q_k^2
\end{equation*}
is almost surely convergent. On the other hand, by a well-known Abel's
result, the convergence of a series $\sum_k a_k$, with $a_k\in\R$,
implies the convergence of the series $\sum_k k^{-1} a_k$ and the
relation $n\sum_{k\geq n} k^{-1}a_k\to 0$ for $n\to +\infty$. Applying
this result, we find (\ref{eq6}) and the proof is so concluded.  \qed

{\bf Proof of Corollary \ref{cor-clt1}. }
  It will suffice to verify
that condition (i) and (ii) of Theorem \ref{clt1} hold with $H=h$.
With regard to condition (ii), it is enough to observe that, 
{by the obvious
inequality $Q_k=r_{k+1}\rho_k\leq c_{k+1}\rho_k$ and the
identity in distribution of the random variables $\rho_k$, we have
\begin{equation*}
\textstyle\sum_{k\geq 0}k^2\,{\rm E}[Q_k^4]
\leq \textstyle\sum_{k\geq 0}k^{2}c_{k+1}^4\,{\rm E}[\rho_k^4]
= {\rm E}[\rho_0^4]\sum_{k\geq 0}k^{2}c_{k+1}^4
<\infty.
\end{equation*} 
}
In order to prove condition (i) of Theorem \ref{clt} (with
$H=h$), we observe that the series
\begin{equation*}
\textstyle\sum_{k}k^{-1}\,(\rho_{k}^2-\beta)
\end{equation*}
is almost surely convergent: indeed, the random variables
$Z_k:=k^{-1}\,(\rho_{k}^2-\beta)$ are independent, centered and
square-integrable, with ${\rm Var}[Z_k]=k^{-2}\,{\rm Var}[\rho_1^2]$.
Therefore, by the above mentioned Abel's result, we obtain the almost
sure convergence of the series
\begin{equation*}
\textstyle\sum_{k}k^{-2}\,(\rho_{k}^2-\beta)
\end{equation*}
and the relation (for $n\to+\infty$)
\begin{equation*}
n\textstyle\sum_{k\geq n}k^{-2}\,(\rho_{k}^2-\beta)
\stackrel{a.s.}\longrightarrow 0.
\end{equation*}
Since we have $n\textstyle\sum_{k\geq n}k^{-2}\to 1$ for $n\to
+\infty$, the above relation can be rewritten in the form
\begin{equation*}
n\textstyle\sum_{k\geq n}k^{-2}\rho_{k}^2
\stackrel{a.s.}\longrightarrow \beta.
\end{equation*}
Now we observe that
\footnote{ Given two sequences $(a_n),\;(b_n)$ of random variables,
the notation $a_n\stackrel{a.s.}\sim b_n$ means that
 $\frac{a_n}{b_n}\stackrel{a.s.}\to 1$.}
\begin{equation*}
Q_k^2=r_{k+1}^2\rho_k^2\stackrel{a.s.}\sim\alpha^2k^{-2}\rho_{k}^2.
\end{equation*}
Hence, for $n\to+\infty$, we have
\begin{equation*}
n\textstyle\sum_{k\geq n}Q_k^2\stackrel{a.s.}\sim
\alpha^{2}n\textstyle\sum_{k\geq n}k^{-2}\,\rho_{k}^2
\stackrel{a.s.}\longrightarrow \alpha^{2}\beta=h.
\end{equation*}
Condition (i) of Theorem \ref{clt1} (with $H=h$) is thus proved
and the proof is concluded.  \qed

\subsection{Proofs of Section \ref{Sec.5}}

In order to study the asymptotic behavior of $(L_n)_{n \geq 1}$ it
will be useful to introduce the sequence of the increments
\begin{equation*}
U_1:=L_1=1\qquad\hbox{and}\qquad
U_{n}:=L_n-L_{n-1}\quad\hbox{for } n\geq 2.
\end{equation*}
Clearly $(U_n)_{n\geq 1}$ is a sequence of random variables with
values in $\{0,1\}$ such that, for each $n\geq 1$, the random variable
$U_n$ is $\CF^X_n$-measurable and $L_n=\sum_{i=1}^n U_n$.

{\bf Proof of Theorem \ref{slln}.}
Without loss of generality, we can assume
$h_n>0$ for each $n$. Let us set
\begin{equation*}
Z_0:=0\qquad Z_n:=\textstyle\sum_{j=1}^n (U_j-s_{j-1})/h_j.
\end{equation*}
Then $Z=(Z_n)_{n\geq 0}$ is a martingale with respect to the filtration ${\cal
F}=({\cal F}_n)_{n\geq 0}$. Indeed,  by Lemma \ref{lemma-prob-cond}, we have
\begin{equation*}
{\rm E}[Z_{n+1}-Z_n\,|\,{\cal F}_{n}]={\rm E}[U_{n+1}-s_n\,|\,{\cal F}_{n}]
={\rm E}[I_{\{L_{n+1}=L_n+1\}}-s_n\,|\,{\cal F}_{n}]
=0. 
\end{equation*} 
Moreover, we have
\begin{equation*}
{\rm E}[U_{n+1}]=P(L_{n+1}=L_n+1)={\rm E}[s_n]
\end{equation*}
and
\begin{equation*}
\begin{split}
{\rm E}\big[(U_{n+1}-s_n)^2\big]&=
{\rm E}\big[\,{\rm E}[(U_{n+1}-s_n)^2\,|\,{\cal F}_n]\,\big]\\
&=
{\rm E}[(1-s_n)^2s_n+s_n^2(1-s_n)]={\rm E}[s_n(1-s_n)].
\end{split}
\end{equation*}
Therefore we obtain
\begin{equation*}
\begin{split}
{\rm E}[Z_{n}^2]=
\textstyle\sum_{j=1}^n {\rm E}\big[(U_j-s_{j-1})^2\big]/h_j^2
=\textstyle\sum_{j=1}^n {\rm E}\big[s_{j-1}(1-s_{j-1})\big]/h_j^2
\end{split}
\end{equation*}
and so 
$${\textstyle\sup_n}{\rm E}[Z_{n}^2]\leq
\textstyle\sum_{j\geq 1}{\rm E}\big[s_{j-1}(1-s_{j-1})\big]/h_j^2<+\infty.
$$ 
It follows that $(Z_n)_{n \geq 1}$ converges almost surely and, by Kronecker's
lemma, we get
\begin{equation*}
\frac{1}{h_n}(L_n-{\textstyle\sum_{j=1}^n} s_{j-1} )=
\frac{1}{h_n}\textstyle\sum_{j=1}^n (U_j-s_{j-1})
\stackrel{a.s.}\longrightarrow 0.
\end{equation*}
Therefore, since
$\sum_{j=1}^n s_{j-1}/h_n=
\sum_{j=0}^n s_{j}/h_n-s_n/h_n\stackrel{a.s.}\longrightarrow L$, we
obtain $L_n/h_n\stackrel{a.s.}\longrightarrow L$.  \qed

In order to prove Theorem \ref{clt} we need a preliminary lemma.

\begin{lemma}\label{lemma-cond-distr}
   If $(X_n)_{n\geq 1}$ is a GOS with $\mu$ diffuse, then (with the
  previous notation), for each fixed $k$, a version of the conditional
  distribution of $(U_j)_{j\geq k+1}$ given $\CG_{k}$ is the
  kernel $Q_k$ so defined:
\begin{equation*}
Q_k(\omega,\cdot):={\textstyle\bigotimes_{j=k+1}^{\infty}}
{\cal B}\big(1,r_{j-1}(\omega)\big)
\end{equation*} 
where ${\cal B}\big(1,r_{j-1}(\omega)\big)$ denotes the Bernoulli
distribution with parameter $r_{j-1}(\omega)$.
\end{lemma}

{\bf Proof.} It is enough to verify that, for each $n\geq
1$, for each $\epsilon_{k+1},\dots,\epsilon_{k+n}\in\{0,1\}$ and for
each ${\CG}_{k}$-measurable real--valued bounded random variable $Z$, we have
\begin{equation}\label{tesi}
{\rm E}\big[Z I_{\{U_{k+1}=\epsilon_{k+1},\dots,U_{k+n}=\epsilon_{k+n}\}}\big]=
{\rm E}\big[Z\, {\textstyle\prod_{j=k+1}^{k+n}} 
r_{j-1}^{\epsilon_j}(1-r_{j-1})^{1-\epsilon_j}
\big].
\end{equation}
We go on with the proof by induction on $n$. For $n=1$, by Lemma 
\ref{lemma-prob-cond}, we have
\begin{equation*}
{\rm E}\big[ZI_{\{U_{k+1}=\epsilon_{k+1}\}}\big]=
{\rm E}\big[Z{\rm E}[I_{\{U_{k+1}=\epsilon_{k+1}\}}\,|\,\CG_k]\,\big]=
{\rm E}[Z r_k^{\epsilon_{k+1}}(1-r_k)^{1-\epsilon_{k+1}}].
\end{equation*}
Assume that (\ref{tesi}) is true for $n-1$ and let us prove it for
$n$. Let us fix an $\CG_{k}$-measurable real--valued bounded random
variable $Z$. By Lemma \ref{lemma-prob-cond}, we have
\begin{equation*}
\begin{split}
&{\rm E}\big[Z
 I_{\{U_{k+1}=\epsilon_{k+1},\dots,U_{k+n}=\epsilon_{k+n}\}}\big]=
{\rm E}\big[Z
 I_{\{U_{k+1}=\epsilon_{k+1},\dots,U_{k+n-1}=\epsilon_{k+n-1}\}}
{\rm E}[U_{k+n}=\epsilon_{k+n}\,|\,\CG_{k+n-1}]\big]\\
&=
{\rm E}\big[Z
r_{k+n-1}^{\epsilon_{k+n}}(1-r_{k+n-1})^{1-\epsilon_{k+n}}
I_{\{U_{k+1}=\epsilon_{k+1},\dots,U_{k+n-1}=\epsilon_{k+n-1}\}}
\big].
\end{split}
\end{equation*}
We have done because also the random variable 
$Zr_{k+n-1}^{\epsilon_{k+n}}(1-r_{k+n-1})^{1-\epsilon_{k+n}}$ is
$\CG_k$-measurable and (\ref{tesi}) is true for $n-1$.
\qed

{\bf Proof of Theorem \ref{clt}.}  Without loss of generality,
we can assume $h_n>0$ for each $n$.  In order to prove the desidered
$\cal A$-stable convergence, it is enough to prove the ${\cal
  F}^X_{\infty}\vee{\cal F}^Y_{\infty}$-stable convergence of $(T_n)$
to ${\cal N}(0,\sigma^2)$. 
But, in order to prove this last
convergence, since we have ${\cal F}^X_{\infty}\vee{\cal
  F}^Y_{\infty}= \bigvee_k{\cal G}_k$, it suffices to prove that, for
each $k$ and $A$ in ${\cal G}_k$ with $P(A)\neq 0$, the sequence
$(T_n)$ converges in distribution under $P_A$ to the probability
measure $P_A{\cal N}(0,\sigma^2)$. In other words, it is sufficient to
fix $k$ and to verify that $(T_{k+n})_n$ (and so $(T_n)_n$) converges
$\CG_k$-stably to ${\cal N}(0,\sigma^2)$. (Note that the kernel ${\cal
  N}(0,\sigma^2)$ is $\CG_k\vee\CN$-measurable for each fixed $k$.)
To this end, we observe that we have
\begin{equation*}
T_{k+n}=\frac{\textstyle\sum_{j=1}^{k+n}(U_j-r_{j-1})}{\sqrt{h_{k+n}}}
=\frac{\textstyle\sum_{j=1}^{k}(U_j-r_{j-1})}{\sqrt{h_{k+n}}}+
\frac{\textstyle\sum_{j=k+1}^{k+n}(U_j-r_{j-1})}{\sqrt{h_{k+n}}}.
\end{equation*}
Obviously, for $n\to +\infty$, we have
\begin{equation*}
\frac{\textstyle\sum_{j=1}^{k}(U_j-r_{j-1})}{\sqrt{h_{k+n}}}
\stackrel{a.s.}\longrightarrow 0.
\end{equation*}
Therefore we have to prove 
\begin{equation}\label{tesi-conv}
\frac{\textstyle\sum_{j=k+1}^{k+n}(U_j-r_{j-1})}{\sqrt{h_{k+n}}}
\stackrel{\CG_k-\hbox{stably}}\longrightarrow {\cal N}(0,\sigma^2).
\end{equation}
\indent From Lemma \ref{lemma-cond-distr} we
know that a version of the the conditional distribution of
$(U_j)_{j\geq k+1}$ given $\CG_{k}$ is
the kernel $Q_k$ so defined:
\begin{equation*}
Q_k(\omega,\cdot)={\textstyle\bigotimes_{j=k+1}^{\infty}}
{\cal B}\big(1,r_{j-1}(\omega)\big).
\end{equation*} 
On the canonical space $\R^{\N^*}$ let us consider the canonical
projections $(\xi_j)_{j\geq k+1}$. Then, for each $n\geq 1$, a version of
the conditional distribution of
\begin{equation*}
\frac{\textstyle\sum_{j=k+1}^{k+n}(U_j-r_{j-1})}{\sqrt{h_{k+n}}}
\end{equation*}
given $\CG_{k}$ is the kernel $N_{k+n}$ so characterized: for
each $\omega$, the probability measure $N_{k+n}(\omega,\cdot)$ is the
distribution, under the probability measure $Q_k(\omega,\cdot)$, of
the random variable (which is defined on the canonical space)
\begin{equation*}
\frac{\textstyle\sum_{j=k+1}^{k+n}\big(\xi_j-r_{j-1}(\omega)\big)}
{\sqrt{h_{k+n}}}.
\end{equation*}
On the other hand, for almost every $\omega$, under $Q_k(\omega,\cdot)$,
the random variables
\begin{equation*} 
Z_{n,i}:=\frac{\xi_{k+i}-r_{k+i-1}(\omega)}{\sqrt{h_{k+n}}}
\qquad\hbox{for } n\geq 1,\; 1\leq i\leq n
\end{equation*}
form a triangular array which satisfies the assumptions of Theorem
\ref{clt-triangular} in the Appendix. Indeed, we have the row-independence
property and
\begin{equation*}
{\rm E}^{Q_k(\omega,\cdot)}[Z_{n,i}]=0,\qquad 
{\rm E}^{Q_k(\omega,\cdot)}[Z_{n,i}^2]=
\frac{r_{k+i-1}(\omega)\big(1-r_{k+i-1}(\omega)\big)}{h_{k+n}}.
\end{equation*}
Therefore, by assumption, for $n\to +\infty$, we have for almost every
$\omega$,
$$
{\textstyle\sum_{i=1}^n} {\rm E}^{Q_k(\omega,\cdot)}[Z_{n,i}^2]= 
\frac{\textstyle\sum_{i=1}^n
  r_{k+i-1}(\omega)\big(1-r_{k+i-1}(\omega)\big)}{h_{k+n}}=
\sigma^2_{k+n}(\omega)-\frac{h_{k-1}\sigma^2_{k-1}(\omega)}{h_{k+n}}
\longrightarrow\sigma^2(\omega).
$$ 
Moreover, under $Q_k(\omega,\cdot)$, 
we have $Z_n^*:=\sup_i Z_{n,i}\leq 2/\sqrt{h_{k+n}}\longrightarrow
0$. Finally, we observe that, setting
$V_n:=\textstyle\sum_{i=1}^n Z_{n,i}^2$, we have
\begin{equation*}
{\rm E}^{Q_k(\omega,\cdot)}[V_n^2]=
{\rm Var}^{Q_k(\omega,\cdot)}[V_n]+
\left(\sigma^2_{k+n}(\omega)-
\frac{h_{k-1}\sigma^2_{k-1}(\omega)}{h_{k+n}}\right)^2
\end{equation*}
with
\begin{equation*} 
\begin{split}
{\rm Var}^{Q_k(\omega,\cdot)}[V_n]=
{\textstyle\sum_{i=1}^n}{\rm Var}^{Q_k(\omega,\cdot)}[Z_{n,i}^2]
&\leq
{\textstyle\sum_{i=1}^n}{\rm E}^{Q_k(\omega,\cdot)}[Z_{n,i}^4]\\
&\leq 4
\left(\sigma^2_{k+n}(\omega)-\frac{h_{k-1}\sigma^2_{k-1}(\omega)}
{h_{k+n}}\right)
  \frac{1}{h_{k+n}}.
\end{split}
\end{equation*}
Since, for almost every $\omega$, the sequence
$(\sigma^2_n(\omega))_n$ is bounded and $h_n\uparrow +\infty$, it
follows that, for almost every $\omega$, the sequence $(V_n)_n$ is
bounded in $L^2$ under $Q_k(\omega,\cdot)$ and so uniformly
integrable.  Theorem \ref{clt-triangular} assures
that, for almost every $\omega$, the sequence of probability measures
\begin{equation*}
\big(N_{k+n}(\omega,\cdot)\big)_{n\geq 1}
\end{equation*}
weakly converges to the Gaussian distribution ${\cal
  N}\big(0,\sigma^2(\omega)\big)$. This fact implies that, for each
  bounded continuous function $g$, we have 
$$
{\rm E}\left[ 
g\left(
\frac{\textstyle\sum_{j=k+1}^{k+n}(U_j-r_{j-1})}{\sqrt{h_{k+n}}}
\right)\,\big|\,
{\CG}_k
\right]
\stackrel{a.s.}\longrightarrow 
{\cal N}(0,\sigma^2)(g).
$$
It obviously follows the $\CG_{k}$-stable convergence
(\ref{tesi-conv}).  \qed

\appendix
\section{Appendix}

For the reader's convenience, we state some results used above.

\begin{theorem}\label{clt-triangular}
Let $(Z_{n,i})_{n\geq 1,\,1\leq i\leq k_n}$ be a triangular array of
square integrable centered random variables on a probability space
$(\Omega,{\cal A},P)$. Suppose that, for each fixed $n$, $(Z_{n,i})_i$
is independent (``row-independence property''). Moreover, set
\begin{equation*}
\begin{split}
&\sigma_{n,i}^2:={\rm E}[Z_{n,i}^2]={\rm Var}[Z_{n,i}],
\qquad\sigma_n^2:=\textstyle\sum_{i=1}^{k_n}\sigma_{n,i}^2,
\\
&V_n:=\textstyle\sum_{i=1}^{k_n} Z_{n,i}^2,
\qquad
Z_n^*:=\textstyle\sup_{1\leq i\leq k_n}|Z_{n,i}|
\end{split}
\end{equation*}
and assume that $(V_n)_{n\geq 1}$ is uniformly integrable,
$Z_n^*\stackrel{P}\longrightarrow 0$ and $\sigma_n^2\longrightarrow
\sigma^2$.  \\[5pt] \indent Then $\sum_{i=1}^{k_n}
Z_{n,i}\stackrel{\hbox{in law}}\longrightarrow{\cal N}(0,\sigma^2)$.
\end{theorem} 
{\bf Proof.} 
 In \cite{hall-heyde} (see pp.
  53--54) it is proved that, under the uniform integrability of
  $(V_n)$, the convergence in probability to zero of $(Z_n^*)_{n\geq 1}$ is
  equivalent to the Lindeberg condition. Hence, it is possible to
  apply Corollary 3.1 (pp. 58-59) in \cite{hall-heyde} with ${\cal
    F}_{n,i}=\sigma(Z_{n,1},\dots,Z_{n,i})$.  

\begin{theorem} \label{fam-tri}
(See Th.~5 and Cor.~7
of sec.~7 in \cite{cri-le-pra})\\
{   Let $(l_n)_{n\geq 1}$ be a sequence of strictly positive
integers. On a probability space $(\Omega,{\cal A},P)$, for each
$n\geq 1$, let $({\cal F}_{n,j})_{0\leq j\leq l_n}$ be a filtration
and $(L_{n,j})_{n\geq 1, 0\leq j\leq l_n}$ be a triangular array of
real random variables such
that, for each $n$, the family $(L_{n,j})_{0\leq j\leq l_n}$ is a
 martingale with respect to $({\cal F}_{n,j})_{0\leq
j\leq l_n}$ and $L_{n,0}=0$.  For each pair $(n,j)$, with $n\geq
1,\,1\leq j\leq l_n$, let us set $Z_{n,j}=L_{n,j}-L_{n,j-1}$ and
\begin{eqnarray*}
S_n=\textstyle\sum_{j=1}^{l_n}Z_{n,j}=L_{n,l_n},
\qquad U_n=\textstyle\sum_{j=1}^{l_n}Z_{n,j}^2,
\qquad Z_n^*=\textstyle\sup_{1\leq j\leq l_n}\;|Z_{n,j}|.
\end{eqnarray*}
Let us suppose that the sequence $(U_n)_{n\geq 1}$ converges in
probability to a positive random variable $U$ and the sequence
$(Z_n^*)_{n\geq 1}$ converges in $L^1$ to zero.  Finally, let $\cal N$
be the sub-$\sigma$-field generated by the $P$-negligible events of
$\cal A$ and let us set
\begin{eqnarray*}
{\cal H}_j=\textstyle\liminf_n{\cal F}_{n,j\wedge l_n}
\quad\hbox{for}\; j\geq 0,
\qquad{\cal H}=\textstyle\bigvee_{j\geq 0}{\cal H}_j.
\end{eqnarray*}
If $U$ is measurable with respect to the $\sigma$-field ${\cal H}\vee{\cal N}$,
then $(S_n)_{n\geq 1}$ converges $\cal H$-stably to the Gaussian
kernel ${\cal N}(0,U)$.  
}
\end{theorem}

\begin{theorem}\label{th-appendice}(see Crimaldi, 2007)\\
{  On~$(\Omega,{\cal A},P)$, let
$(V_n)_{n\geq 0}$ be a real martingale with respect to a filtration
$\CG=(\CG_n)_{n\geq 0}$. Suppose that $(V_n)_{n\geq 0}$ converges in
$L^1$ to a random variable $V$. Moreover, setting
\begin{equation}\label{eq-def}
\textstyle U_n:=n\sum_{k\geq n}(V_k-V_{k+1})^2,\qquad
\textstyle Z:=\sup_k\sqrt{k}\,|V_k-V_{k+1}|,
\end{equation}
assume that the following conditions hold:
\\[5pt]
\indent (a) The random variable $Z$ is integrable.
\\[5pt]
\indent (b) The sequence $(U_n)_{n\geq 1}$
converges almost surely to a positive real random variable~$U$.
\\[5pt] \indent Then, with respect to $\CG$, the sequence $(W_n)_{n\geq 1}$ 
defined by
\begin{equation}\label{eq-def2}
\textstyle W_n:=\sqrt{n}(V_n-V)
\end{equation}
converges to the Gaussian kernel~${\cal N}(0,U)$ in the sense of the
almost sure conditional convergence.}
\end{theorem}

  Obviously the previous almost sure conditional convergence also
  holds with respect to any  filtration~$\CF'$ such that
  $\CF^V_n\subset\CF'_n\subset\CG_n$.

\end{document}